\documentclass[authoryear,webpdf,imanum]{ima-authoring-template}%
\usepackage{amsthm}
\graphicspath{{Fig/}}

\theoremstyle{thmstyletwo}%
\newtheorem{theorem}{Theorem}[section]
\usepackage{bm}
\newtheorem{remark}{Remark}%

\numberwithin{equation}{section}
\newtheorem{lemma}[theorem]{Lemma}

\newcommand{\ii}{\mathrm i}
\begin{document}

\copyrightyear{2026}
\vol{00}
\pubyear{}
\access{}
\appnotes{Paper}
\copyrightstatement{Published by Oxford University Press on behalf of the Institute of Mathematics and its Applications. All rights reserved.}
\firstpage{1}


\title[An explicit third-order low-regularity integrator]{An explicit third-order low-regularity integrator for the cubic nonlinear Schr\"odinger equation}

\author{Shengxin Luo
\address{\orgdiv{School of Mathematics and Statistics}, \orgname{Xi'an Jiaotong University}, \postcode{710049}, \state{Xi'an}, \country{China}}}
\author{Bin Wang*
\address{\orgdiv{School of Mathematics and Statistics}, \orgname{Xi'an Jiaotong University}, \postcode{710049}, \state{Xi'an}, \country{China}}}
\author{Yan Wang
\address{\orgdiv{School of Mathematics and Statistics}, \orgname{Central China Normal University}, \postcode{430079}, \state{Wuhan}, \country{China}}}

\authormark{Luo et al.}

\corresp[*]{Corresponding author: \href{wangbinmaths@xjtu.edu.cn}{wangbinmaths@xjtu.edu.cn}}

\received{Date}{0}{Year}
\revised{Date}{0}{Year}
\accepted{Date}{0}{Year}


\abstract{In this paper, we propose an explicit third-order
	low-regularity integrator for the cubic nonlinear Schr\"odinger
	equation on the one-dimensional torus $\mathbb{T}$.
	By carefully designing numerical approximations for the integrals
	in Fourier space and exploiting oscillatory cancellations across
	time steps, we establish rigorous global error estimates under
	low regularity assumptions.
	More precisely, we prove third-order convergence in
	$H^r(\mathbb{T})$ for $r>1/2$ with initial data in
	$H^{r+2}(\mathbb{T})$.
	Furthermore, third-order convergence in $L^2(\mathbb{T})$
	is established for initial data in $H^{2}(\mathbb{T})$.
	Numerical experiments are presented to confirm the theoretical
	convergence results.}


\maketitle

\section{Introduction}\label{sec1}

As a prototypical dispersive equation, the nonlinear Schr\"odinger (NLS) equation has found extensive applications in nonlinear optics, plasma physics, and Bose–Einstein condensation  \cite{sulem1999,agrawal2019nonlinear}. In this paper, we are concerned with the numerical approximation of the cubic nonlinear Schr\"odinger equation
\begin{equation}\label{NLS}
	\begin{aligned}
		&\ii\,\partial_t u(t,x)+\partial_x^{2}u(t,x)
		= \lambda\,|u(t,x)|^{2}u(t,x),\ \ u(0,x) = u^{0}(x), \ \  x\in\mathbb{T},\ t\in(0,T],
	\end{aligned}
\end{equation}
on the torus $\mathbb{T}=[-\pi,\pi]$ with periodic boundary condition. The initial data satisfies 
$u^{0}(x)\in H^{s}(\mathbb{T})$ for $ s\ge 0$, and the  values $\lambda=-1$ and $1$ correspond to the focusing and defocusing case of NLS, respectively. The local well-posedness of the problem \eqref{NLS} in $H^s(\mathbb{T})$ is well established for $s \ge 0$ \cite{bourgain1993fourier,Ostermann2022}. 
\par
Numerical schemes for the NLS equation have been extensively studied. Standard approaches include operator splitting methods (e.g., Lie~\cite{ignat2011splitting,su2020time} and Strang splitting~\cite{eilinghoff2016fractional,lubich2008splitting}), exponential integrators~\cite{berland2007conservation,celledoni2008symmetric,hochbruck2010exponential}, and finite difference schemes~\cite{akrivis1993finite,henning2017crank,sanzserna1984methods}. These classical methods typically achieve their expected convergence orders provided that the solution of \eqref{NLS} is sufficiently smooth. For instance, 
to achieve second-order convergence in $H^s(\mathbb{T})$ ($s\ge0$) with the Strang splitting scheme, the solution is required to be in $ C([0,T];H^{s+4}(\mathbb{T}))$ (see \cite{lubich2008splitting}). Generally, for a classical scheme to achieve $k$th-order convergence in $H^{s}(\mathbb{T})$ ($s\ge 0$), the regularity requirement for the solution is $u\in C([0,T];H^{s+2k}(\mathbb{T}))$, as well as for the initial data $u^{0}(x)$. 

However, in many practical applications, the initial data are naturally rough and do not satisfy the high-regularity assumptions underlying classical convergence analyses.
A notable example in ocean engineering is the occurrence of freak waves in irregular random sea states, posing serious hazards to ships and offshore structures~\cite{Onorato2001}.
In this setting, a random sea state characterized by a JONSWAP power spectrum is considered and, assuming the dynamics to cubic order, the leading-order evolution is modeled by the NLS equation.
The initial conditions are generated via a random-phase spectral representation, with phases taken as independent uniform random variables on $(0,2\pi)$ \cite{Onorato2001}.
Such a random-phase initialization produces intrinsically irregular initial wave fields, placing the problem in a low-regularity regime where standard time-stepping schemes often exhibit order reduction and classical smoothness-based analyses may no longer apply.
Consequently, reliable numerical simulations of the NLS equation with rough initial data become essential.

In recent years, various low-regularity integrators have been developed for the NLS equation to alleviate the derivative-loss in the consistency error and  maintain reliable accuracy under weaker smoothness assumptions on the solution. The pioneering scheme proposed by~\cite{ostermann2018low} is constructed based on twisted variable (i.e., $v(t) = \mathrm e^{-\ii t\partial_x^2} u(t)$) and harmonic analysis techniques. It was proven to achieve first-order convergence in $H^s(\mathbb{T}^d)$ under the assumption $u^0 \in H^{s+1}(\mathbb{T}^d)$. Subsequently, different harmonic analysis techniques were employed in \cite{li2021fully} to construct a fully discrete low-regularity integrator for the one-dimensional periodic cubic NLS equation, establishing an almost first-order convergence rate in $L^{2}(\mathbb{T})$ for initial data in $H^{1}(\mathbb{T})$. Later, another low-regularity integrator was introduced in \cite{Wu2021}, which achieves first-order convergence in $H^{s}(\mathbb{T})$ without any loss of regularity. More recently, the fully discrete low-regularity analysis was refined in \cite{Ostermannyao2022} by allowing initial data $u^{0}\in H^{s}(\mathbb{T})$ for $s\in(1/2,1]$ while preserving the almost first-order convergence at the endpoint $s=1$. Specifically, the integrator in \cite{Ostermannyao2022} attains a convergence rate of order $\tfrac{3}{2}s-\tfrac{1}{2}-\varepsilon$ in $L^{2}(\mathbb{T})$ for $s\in(1/2,1]$ and a small quantity $\varepsilon>0$, lowering the regularity threshold required in \cite{li2021fully}.


Recent work has also made significant progress in developing second-order schemes for rough initial data. 
A low-regularity integrator for the cubic NLS equation was proposed in \cite{ostermann2022second}, establishing second-order convergence in $H^s(\mathbb{T}^d)$ provided that $u^0 \in H^{s+2}(\mathbb{T}^d)$. 
More recently, \cite{cao2024new} further reduced the regularity requirement to $u^0 \in H^{5/3}(\mathbb{T})$ for second-order convergence in $L^2(\mathbb{T})$.
This analytical framework, built around a twisted variable reformulation and supported by harmonic analysis estimates, has also been adapted to other nonlinear dispersive equations. 
For instance, related low-regularity schemes have been developed for the Korteweg--de Vries, Dirac, Klein--Gordon, and Boussinesq equations  (see, e.g., \cite{Hofmanov2017,Ning2022,Ostermann2020,Ostermannsu2019,Li2023,Schratz2021,Wang2022,Wuzh2022,Wuzh22022}).  
Beyond differential equations with rough initial values, numerical algorithms for other low-regularity problems, such as PDEs with low-regularity potentials, have been intensively studied in recent years \cite{Zhao_2021,Bao_2024,Bao_2024_2}.  

\par Beyond second order, the theoretical possibility of constructing low-regularity integrators of arbitrarily high order has been established in \cite{bruned2022resonance}.
Nevertheless, translating these general construction principles into explicit higher order schemes in the low-regularity regime remains nontrivial, as it necessitates simultaneously overcoming algebraic complexity and regularity restrictions. Specifically, this task faces two primary obstacles.
First, achieving third-order convergence often requires tracking many coupled correction terms.
To obtain a local fourth order error bound for an oscillatory Duhamel integral, one must account for the variation of the nonlinearity over the time step, rather than freezing it at some fixed $t$.
This yields many interacting contributions, where the $\varphi$-functions couple nontrivially with nonlinear evaluations and intermediate stages, making the derivation and the subsequent stability and local error analysis more involved.
Second, there is an intrinsic tension between achieving third-order convergence and low regularity.
Standard Taylor expansions of oscillatory factors in Fourier space, such as $\mathrm e^{\ii\tau k^{2}}$ with the time step size $\tau$, inevitably generate powers of $k^{2}$, which correspond to the spatial derivative $\partial_x^2$.
These terms amplify high frequencies and are not uniformly bounded for low-regularity data, leading to derivative loss, order reduction, and the breakdown of low-regularity convergence.
\par 
%
To overcome these difficulties, we employ three main techniques:
\begin{itemize}
	\item Substitute the second-order low-regularity approximation into
	the Duhamel integral for the twisted variable $v$.
	\item Decompose the oscillatory phases algebraically to distribute
	spatial derivatives more evenly, and introduce Fourier filtering
	to mitigate derivative loss.
	\item Exploit discrete oscillatory cancellation and resonance counting
	to control the partial sums of the cubic and quintic remainders
	for fixed inputs, and apply discrete summation by parts to handle
	time-dependent inputs.
\end{itemize}
Building on these strategies,  we develop an explicit third-order low-regularity integrator (denoted by \texttt{LRI3}) for the cubic NLS equation and establish its convergence under reduced smoothness assumptions. In particular, we prove global third-order convergence in
$H^r(\mathbb{T})$ for $r>1/2$ assuming
$u\in C([0,T];H^{r+2}(\mathbb{T}))$, and global third-order
convergence in $L^2(\mathbb{T})$ with the requirement
$u\in C([0,T];H^{2}(\mathbb{T}))$.

\par The remainder of this paper is organized as follows.
Section~\ref{preliminary} introduces the notation and reviews
the necessary preliminaries.
Section~\ref{res} presents the main results, including the
\texttt{LRI3} scheme, its convergence theorem, and numerical experiments.
Section~\ref{construction} describes the derivation of \texttt{LRI3}
and establishes the key remainder estimates underlying the convergence analysis.
Section~\ref{proof} provides the proofs of the technical lemmas
and the main convergence theorem.
Finally, Section~\ref{section:conclusion} presents some concluding remarks.
\par
\section{Preliminaries}
\label{preliminary}
In this section, we introduce the notation, recall some basic facts
from harmonic analysis, and state several auxiliary lemmas used
in the subsequent analysis.
\subsection{Some notations}
We write $A \lesssim B$  to denote $A \le CB$ for a positive constant $C$ which may vary from line to line but is always independent of the time step size $\tau>0$ and $n=0,1,\ldots, T/\tau$. 
For $f\in L^{2}(\mathbb{T})$, we write
$f(x)=\sum_{k\in\mathbb{Z}}\hat f_k\mathrm e^{\ii kx}$ with
$\hat f_k=\mathcal{F}_k[f]:=\frac{1}{2\pi}\int_{\mathbb{T}}\mathrm e^{-\ii kx}f(x)\,\mathrm{d}x$, $k\in\mathbb{Z}$.
The following properties of the Fourier transform are used frequently:
\begin{align*}
\|f\|_{L^{2}}^{2}
= 2\pi \sum_{k\in\mathbb{Z}} |\hat{f}_{k}|^{2}
\ \ \text{(Plancherel identity),}\ \
\mathcal{F}_{k}[fg]
= \sum_{k_{1}\in\mathbb{Z}} \hat{f}_{k-k_{1}}\,\hat{g}_{k_{1}}
\ \  \text{(Convolution).}
\end{align*}
Consistent with this notation, the zero-mode and nonzero-mode operators  are respectively defined by
$
\Pi_{0} f = \hat{f}_{0}$ and $
\Pi_{\neq 0} f = \sum_{k \in \mathbb{Z},k \neq 0} \hat{f}_{k}\,\mathrm e^{\ii kx} .$ 
For notational convenience, we denote
$
\langle k\rangle=(1+k^{2})^{\frac{1}{2}}$ and $
J=\langle \ii^{-1}\partial_{x}\rangle.
$
With this notation, an equivalent norm on $H^{s}$ is given by
$$
\|f\|_{H^{s}}^{2}:=\|J^{s}f\|_{L^{2}}^{2}
=2\pi\sum_{k\in\mathbb{Z}}(1+k^{2})^{s}\,|\hat f_{k}|^{2}.$$
\par Using the Fourier expansion, for any $t\in\mathbb{R}$, the action of the free Schr\"odinger group $\mathrm e^{\ii t\partial_x^{2}}$ can be written as
$
\mathrm e^{\ii t\partial_{x}^{2}}f(x)
=\sum_{k\in\mathbb{Z}}\mathrm e^{-\ii t k^{2}}\,
\hat f_{k}\,\mathrm e^{\ii k x}.$
Since $|\mathrm e^{-\ii t k^{2}}|=1$, $\mathrm e^{\ii t\partial_{x}^{2}}$ is a linear isometry on $H^{s}(\mathbb{T})$ for $s\ge0$. 
Moreover, for any integer $m\ge1$ we define a regularization of
$\partial_{x}^{-m}:H^{s}\to H^{s+m}$ by $\mathcal{F}_{k}[\partial_{x}^{-m}f]=(\ii k)^{-m}\,\hat f_{k}$ for $k\neq 0$
and $\mathcal{F}_{k}[\partial_{x}^{-m}f]=0$ for $k=0$. 
We also introduce the following families of functions
for the definition of the numerical scheme:
\begin{equation}\label{adefin}
	\begin{gathered}
		\varphi_j(z)
		=\frac{1}{(j-1)!}\int_0^1
		(1-s)^{j-1}\mathrm e^{sz}\,\mathrm{d}s,
		\ \ \psi_j(z)
		=\frac{(-1)^j}{j!}\int_0^1
		s^j \mathrm e^{sz}\,\mathrm{d}s,
		\ \ j=1,2,3,\\
		\Phi_j(z)=\int_0^1(1-s)
		\big(1-\frac{j}{2}(1-s)\big)
		\mathrm e^{sz}\,\mathrm{d}s,
		\ \ j=1,2.
	\end{gathered}
\end{equation}
These functions are bounded on the imaginary axis, and hence
the corresponding Fourier multipliers are uniformly bounded
on $H^s(\mathbb{T})$ for $s\ge0$ and $t\in\mathbb{R}$.
Throughout, we abbreviate $H^s(\mathbb{T})$ and
$L^2(\mathbb{T})$ as $H^s$ and $L^2$, respectively.
\subsection{Some estimates}
We use the following form of the Kato–Ponce inequality proved in \cite{kato1988commutator}.
\begin{lemma}\label{kp}
	\leavevmode\par
\begin{enumerate}[label=(\roman*)]
\item
For any $\gamma>\frac{1}{2}$ and $f,g\in H^{\gamma}$, we have
$
\|J^{\gamma}(fg)\|_{L^{2}} \lesssim\|f\|_{H^{\gamma}}\|g\|_{H^{\gamma}}.$
\item
For any $\delta \ge 0$, $\gamma>\frac{1}{2}$ and $f\in H^{\delta+\gamma}$, $g\in H^{\delta}$, it is true that
$
\|J^{\delta}(fg)\|_{L^{2}} \lesssim \|f\|_{H^{\delta+\gamma}}\|g\|_{H^{\delta}}.$
\end{enumerate}
\end{lemma}

The following lemma is useful in our analysis. We refer to
\cite{cao2024new} for its proof.
\begin{lemma}\label{li}
	\leavevmode\par
\begin{enumerate}[label=(\roman*)]
\item
If $s>\tfrac{1}{2}$ and $f,g\in H^{s}$, one has
$
\bigl\| J^{-1}{((J f)\cdot g)} \bigr\|_{H^{s}}
\lesssim
\|f\|_{H^{s}}\|g\|_{H^{s}}.$
\item
If $f,g\in H^{1}$, it is obtained that
$
\bigl\| J^{-1}((J f)\cdot g) \bigr\|_{L^{2}}
\lesssim
\min\left\{ \|f\|_{L^{2}}\|g\|_{H^{1}},
\|g\|_{L^{2}}\|f\|_{H^{1}} \right\}.$
\end{enumerate}
\end{lemma}

The following lemma will be used later. Its proof is straightforward and therefore omitted.
\begin{lemma}\label{jb}
For any $\alpha \in \mathbb{R}$, it holds that
$
|\mathrm e^{\ii\alpha}-1|
\lesssim \min\{1,|\alpha|\}
\le \frac{2|\alpha|}{1+|\alpha|}$ and $|\mathrm e^{\ii\alpha}-1-\ii\alpha|
\lesssim\min\{|\alpha|^2,|\alpha|\}\le \frac{2|\alpha|^2}{1+|\alpha|}.$ 
\end{lemma}
\begin{lemma}\label{lem:discrete-time-summation}
	Let $t_m=m\tau$. Define
	\[
	g_\tau(\phi)=\min\left\{1,\frac{\tau}{|1-\mathrm e^{\ii\tau\phi}|}\right\},
	\]
	with $g_\tau(\phi)=1$ when $e^{\ii\tau\phi}=1$. Then the following
	statements hold.
	\begin{enumerate}
		\item For $0\le n\le L=T/\tau$,
		\[
		\left|\sum_{m=0}^{n-1}\mathrm e^{\ii t_m\phi}\right|
		\lesssim_T  \tau^{-1}g_\tau(\phi).
		\]
		
		\item Let
		\[
		I_\ell=\left[\frac{2\pi\ell-\pi}{\tau},
		\frac{2\pi\ell+\pi}{\tau}\right),
		\qquad \ell\in\mathbb Z.
		\]
		If $P$ is a real quadratic polynomial whose leading coefficient
		has modulus at least one, then
		\[
		\sum_{\substack{j\in\mathbb Z\\P(j)\in I_\ell}}
		g_\tau(P(j))^\rho\lesssim_\rho 1
		\qquad \text{for } \rho>\tfrac12.
		\]
		If $P$ is affine with slope of modulus at least one, the same
		bound holds for $\rho>1$. The constants are uniform in $\ell$,
		$\tau$, and the coefficients of $P$ under these assumptions.
		
		\item Let $\mathcal D_m:X^d\to Y$ be real multilinear maps
		between Banach spaces, and set
		$\mathcal A_n=\sum_{m=0}^{n-1}\mathcal D_m$. Suppose that
		\[
		\max_{0\le n\le L}
		\|\mathcal A_n(f_1,\ldots,f_d)\|_Y
		\le M_\tau\prod_{j=1}^d\|f_j\|_X.
		\]
		For $w_0,\ldots,w_{L-1}\in X$, let
		$R=\max\limits_{0\le m <L}\|w_m\|_X$. Then, for $1\le n\le L$,
		\[
		\left\|\sum_{m=0}^{n-1}\mathcal D_m(w_m,\ldots,w_m)\right\|_Y
		\le M_\tau\left(
		R^d+dR^{d-1}\sum_{m=0}^{n-2}\|w_{m+1}-w_m\|_X
		\right).
		\]
	\end{enumerate}
\end{lemma}

\begin{proof}
	The geometric sum is bounded by $n\le L$ and, when
	$\mathrm e^{\ii\tau\phi}\ne1$, by $2/|1-\mathrm e^{\ii\tau\phi}|$. This gives
	the first statement.
	
	For $\phi\in I_\ell$, let $h=|\phi-2\pi\ell/\tau|$. Since
	$\tau h\le\pi$, we have $g_\tau(\phi)\lesssim(1+h)^{-1}$.
	Writing $P(x)=a(x-x^*)^2+b$, we see that on each branch the
	preimage of a phase interval of length $D$ is longest when
	the interval has the vertex value as an endpoint and lies
	in the range of $P$. Its length is therefore at most
	$\sqrt{D/|a|}\le\sqrt D$, and the two branches together contain
	at most $2(1+\sqrt D)$ integers. For an affine $P$ with slope
	of modulus at least one, the preimage has length at most $D$
	and hence contains at most $1+D$ integers.
	Applying these counts to $h\le1$ and to the two phase intervals
	in each layer $2^q<h\le2^{q+1}$ yields
	\[
	\sum_{\substack{j\in\mathbb Z\\P(j)\in I_\ell}}
	g_\tau(P(j))^\rho
	\lesssim_\rho
	1+\sum_{q\ge0}2^{-q\rho}(1+2^{q\theta}),
	\]
	where $\theta=\tfrac12$ in the quadratic case and $\theta=1$
	in the affine case. The series converges for $\rho>\theta$.
	
	Finally, write $\boldsymbol w_m=(w_m,\ldots,w_m)$. Since
	$\mathcal D_m=\mathcal A_{m+1}-\mathcal A_m$, summation by parts
	gives
	\[
	\sum_{m=0}^{n-1}\mathcal D_m(\boldsymbol w_m)
	=\mathcal A_n(\boldsymbol w_{n-1})
	+\sum_{m=0}^{n-2}
	\bigl[\mathcal A_{m+1}(\boldsymbol w_m)
	-\mathcal A_{m+1}(\boldsymbol w_{m+1})\bigr].
	\]
	Expanding each difference one input at a time gives $d$ terms,
	each with one factor $w_m-w_{m+1}$ and $d-1$ factors bounded
	by $R$ in $X$. The assumed bound for $\mathcal A_n$ proves
	the last statement.
\end{proof}

\section{Main results}\label{res}
In this section, we introduce the novel third-order low-regularity integrator \texttt{LRI3}, state the convergence result, and present a numerical experiment validating the theoretical findings.
\subsection{A new LRI3 and its convergence}
Let $t_n = n\tau$, $n = 0, 1, \ldots, L$, be a uniform partition of the time interval $[0,T]$ with step size 
$\tau = T/L$, where $L$ is a positive integer. We define the explicit third-order low-regularity integrator  (\texttt{LRI3}) as
\begin{equation}\label{littlescheme}
u^{n+1} =\widetilde{\Psi}_{\tau}(u^n), \ \ n = 0,1,\ldots,L-1,
\end{equation}
where 
\begin{align}\label{numericalscheme}
\widetilde{\Psi}_{\tau}(f)= \tau^{-1}\widetilde{\Psi}^{[-1]}_{\tau}(f)+ \widetilde{\Psi}^{[0]}_{\tau}(f)+ \tau\widetilde{\Psi}^{[1]}_{\tau}(f)+\tau^2\widetilde{\Psi}^{[2]}_{\tau}(f)+\tau^3\widetilde{\Psi}^{[3]}_{\tau}(f).
\end{align}
For notational simplicity, we set $S_t=\mathrm e^{\ii t\partial_x^2}$,
$D=\partial_x$, and $f_\tau=P_\tau f$. For $a=1,2$, we write $\varphi_{j,a}^{\pm}=\varphi_j(\pm a\ii\tau\partial_x^2)$ for $j=1,2,3$ and $\psi_{j,a}^{\pm}=\psi_j(\pm a\ii\tau\partial_x^2)$ for $j=1,2$. We also set $\Phi_{j,2}^-=\Phi_j(-2\ii\tau\partial_x^2)$ for $j=1,2$.
With this notation, the operators \(\widetilde{\Psi}^{[j]}_\tau\), \(j=-1,\ldots,3\), are given by
\begingroup
\allowdisplaybreaks[4]
\begin{align*}
	\widetilde{\Psi}^{[-1]}_{\tau}(f)&:=\frac{\ii\lambda}{4}D^{-3}\Bigl\{ S_{-\tau}D^{-1}\bar{f}\cdot
	(S_\tau f)^2 - S_{-\tau}D^{-1}\bar{f}\cdot S_\tau f^2 \Bigr\} - \frac{\ii\lambda}{4} S_\tau D^{-3}\Bigl\{ D^{-1}\bar{f}\cdot
	S_{-\tau}(S_\tau f)^2 - D^{-1}\bar{f}\cdot f^2 \Bigr\},\\
	\widetilde{\Psi}^{[0]}_\tau(f)&:=\frac{\lambda}{2}\Pi_0\Bigl\{ \bar{f}\cdot S_{-\tau}(S_\tau D^{-1} f)^2
	- \bar{f}(D^{-1}f)^2 \Bigr\}+\lambda D^{-1}\Bigl\{S_\tau f \cdot D^{-1}\!\big(|S_\tau f|^2\big)
	- S_\tau\big[f\cdot D^{-1}(|f|^2)\big]\Bigr\} \notag \\
	&\quad - \frac{\lambda}{2}D^{-2}\Bigl\{ S_{-\tau}\bar{f}\cdot S_\tau f^2
	- S_\tau \big(f\cdot |f|^2\big)\Bigr\} - \frac{\lambda}{2} S_\tau D^{-1}
	\Bigl\{D \bar{f} \cdot S_{-\tau}\big(S_\tau D^{-1} f\big)^2
	- D \bar f \cdot (D^{-1}f)^2 \Bigr\} \\
	&\quad - \frac{\lambda}{2}D^{-2}\Bigl\{ \Pi_{\neq 0}(S_{-\tau}\bar{f}) \cdot (S_\tau f)^2
	- \Pi_{\neq 0}(S_{-\tau}\bar{f}) \cdot S_\tau f^2 \Bigr\}+S_\tau f,\\
	\widetilde{\Psi}^{[1]}_\tau(f)&:=\ii\lambda |\Pi_0(f)|^2\Pi_0(f)-2\ii\lambda\Big\{ \Pi_0(f D \bar{f})D^{-1}
	+\Pi_0(|f|^2) \Big\} S_\tau f - \ii\lambda (\Pi_0(f))^2 \Pi_{\neq 0}\big(S_\tau \bar{f}\big)
	 \notag \\
	&\quad + 2\ii\lambda \Pi_0(f)\cdot S_\tau D^{-1}(D \bar f \cdot f)- \frac{\ii\lambda}{2} S_\tau D^{-1}\Bigl\{ D \bar f \cdot
S_{-\tau}(S_\tau f)^2 - D\bar f \cdot f^2 \Bigr\}  - \ii\lambda S_\tau D^{-1}(D \bar f \cdot f^2),\\
	\widetilde{\Psi}^{[2]}_\tau(f)&:=2\lambda^2 S_\tau\Big\{
	f^2|f|^2
	\Big[\big(\psi_{1,2}^-+\tfrac12-\varphi_{2,2}^-\big)
	\bar f\Big]\Big\}+\lambda^2 S_\tau\Big\{
	|f|^4(\varphi_{2,2}^-f)-f^3\big(\psi_{1,2}^-+\tfrac12\big)
	\{\bar f\,^2\}\Big\},\\
	\widetilde{\Psi}^{[3]}_\tau(f)&:=4\ii\lambda^3 S_\tau \Bigl\{\big(\psi_{2,2}^-\overline{f}\big) |f|^4 f^2-\psi_{2,2}^-\bigl\{|f|^2 \overline{f}\bigr\}|f|^2 f^2+\frac{1}{4} \psi_{2,2}^-\bigl\{|f|^4 \overline{f}\bigr\}f^2\Bigr\}\\
	&\quad+2\ii\lambda^2 S_\tau\Big\{
	f_\tau^3\big(\psi_{2,2}^--\tfrac16\big)\big\{(D\bar f_\tau)^2\big\}
	-D\{f_\tau^2\}D\big[
	2f_\tau\big(\psi_{2,2}^--\tfrac16\big)\{\bar f_\tau^{\,2}\}
	+\bar f_\tau^{\,2}(\Phi_{1,2}^-f_\tau)\big]\Big\}\\
	&\quad+2\ii\lambda^2 S_\tau\Big\{f_\tau^2(D\bar f_\tau)^2(\varphi_{3,2}^-f_\tau)-(Df_\tau)^2\big[
	2f_\tau\big(\psi_{2,2}^--\tfrac16\big)\{\bar f_\tau^{\,2}\}
	+\bar f_\tau^{\,2}(\Phi_{1,2}^-f_\tau)\big]\Big\}\\
	&\quad+4\ii\lambda^2 S_\tau\Big\{	2f_\tau|f_\tau|^2(Df_\tau)
	\big[\big(\psi_{2,2}^--\tfrac16+\Phi_{2,2}^-\big)D\bar f_\tau\big]+|f_\tau|^2(Df_\tau)^2
	\big[\big(\psi_{2,2}^--\tfrac16+\Phi_{2,2}^-\big)\bar f_\tau\big]\Big\}\\
	&\quad+4\ii\lambda^2 S_\tau\Big\{
	D(|f_\tau|^2)D\big[f_\tau^2\big(\Phi_{1,2}^-\bar f_\tau\big)\big]	+f_\tau|f_\tau|^2(D\bar f_\tau)
	\big[\big(\varphi_{3,2}^--2\Phi_{1,2}^-\big)Df_\tau\big]\Big\}\\
	&\quad+4\ii\lambda^2 S_\tau\Big\{f_\tau^2(Df_\tau)(D\bar f_\tau)(\Phi_{1,2}^-\bar f_\tau)-3f_\tau^2(Df_\tau)\big(\psi_{2,2}^--\tfrac16\big)
	\big[\bar f_\tau(D\bar f_\tau)\big]\Big\}\\
		&\quad+8\ii\lambda^2 S_\tau\Big\{
	D\big[f_\tau\big((\psi_{2,2}^--\tfrac16)\bar f_\tau\big)\big]
	D(f_\tau|f_\tau|^2)
	+f_\tau^2(Df_\tau)(D\bar f_\tau)
	\big[\big(\psi_{2,2}^--\tfrac16\big)\bar f_\tau\big]\Big\}.
\end{align*}
\endgroup
\par The primary theoretical result of this paper is stated in the following theorem.
\begin{theorem}\label{thm1.1}
Let $u^n$ be the numerical solution of the NLS equation given by \eqref{littlescheme}--\eqref{numericalscheme} and $r>\tfrac{1}{2}$.
Then, under the regularity condition $u^0 \in H^{r+2}$, 
there exist positive constants $\tau_0$ and $C_0$ such that, for any 
$\tau \in (0,\tau_0]$, it holds
\begin{equation}\label{T1.1}
\max_{1 \leq n \leq L} \| u(t_n,\cdot) - u^n \|_{H^r} \leq C_0 \tau^3,
\end{equation}
where $C_0$ and $\tau_0$ depend only on the final time $T$ 
and $\| u^0\|_{H^{r+2}}$.  Moreover, under the regularity condition 
$u^0 \in H^{2}$, there exist positive constants $\tau_1$ and $C_1$ such that, for any $\tau \in (0,\tau_1]$, the following error estimate holds:
\begin{equation}\label{T1.2}
\max_{1 \leq n \leq L} \| u(t_n,\cdot) - u^n \|_{L^2} 
\leq C_1 \tau^3,
\end{equation}
where $\tau_1$ and $C_1$ depend only on $T$ and
$\|u^0\|_{H^{2}}$.
\end{theorem}
\begin{remark}\label{r1:lg}
For standard numerical methods, classical convergence analysis typically dictates that achieving the third-order accuracy \eqref{T1.1} requires $u^{0}\in H^{r+6}$.  The proposed new method \texttt{LRI3},  however, circumvents this  restriction to $u^{0}\in H^{r+2}$. This reduction is significant for the scientific computation of the  NLS equation  with rough solutions.
\end{remark}

\begin{remark}\label{r2:fourier}
For practical implementation issues, we employ the spectral method for the spatial discretization, which ensures spectral accuracy due to the periodicity of the domain. As a result, this paper focuses only on the temporal error analysis. The fully discrete error estimates follow readily by coupling these temporal results with standard spectral error bounds.
\end{remark}

\subsection{Numerical test}
In this subsection, we perform a numerical test to confirm the predicted convergence of \eqref{NLS} for the case $\lambda=1$, comparing the proposed method with the Strang splitting, the symmetric second-order low-regularity integrator \texttt{LRI2} from \cite{bronsard2025symmetric}, and the classical third-order exponential integrator \texttt{EI3} from \cite{hochbruck2005explicit}. 
In our numerical experiments, we use a Fourier spectral method with the largest mode $K=2^{12}$ (spatial mesh $\Delta x = 2\pi/2^{13}$) and construct rough initial data $u^0 \in H^\theta$ in the frequency domain as
$\hat{u}_k^0 = \xi_k / \omega_k$, $k=-K+1,\dots,K$, 
where $\omega_k = |k|^{\theta+0.51}$ for $k\neq 0$, $\omega_0=1$, and $\xi_k$ are independent complex uniform random variables in the unit square. Finally, the initial data are normalized in $L^2$ as $u^0 \mapsto u^0 / \|u^0\|_{L^2}$.
\par To evaluate the accuracy of the proposed scheme, we compute a reference solution using the second-order low-regularity integrator of \cite{bruned2022resonance} with time step $\tau_{\mathrm{ref}}=10^{-5}$, and use it to compare the numerical performance of our third-order integrator with the Strang splitting method, the classical third-order exponential integrator \texttt{EI3}, and the second-order low-regularity integrator \texttt{LRI2}.
\begin{figure}[htbp]
$$\begin{array}{cc}
\includegraphics[width=0.33\textwidth]{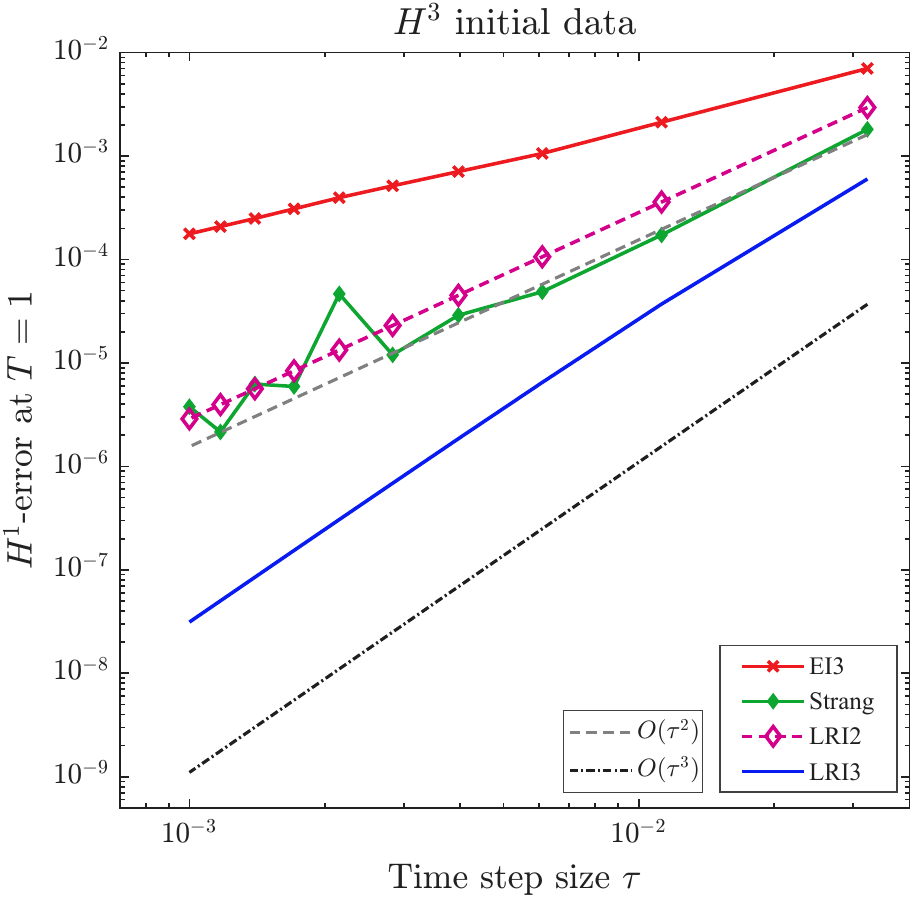}
\includegraphics[width=0.33\textwidth]{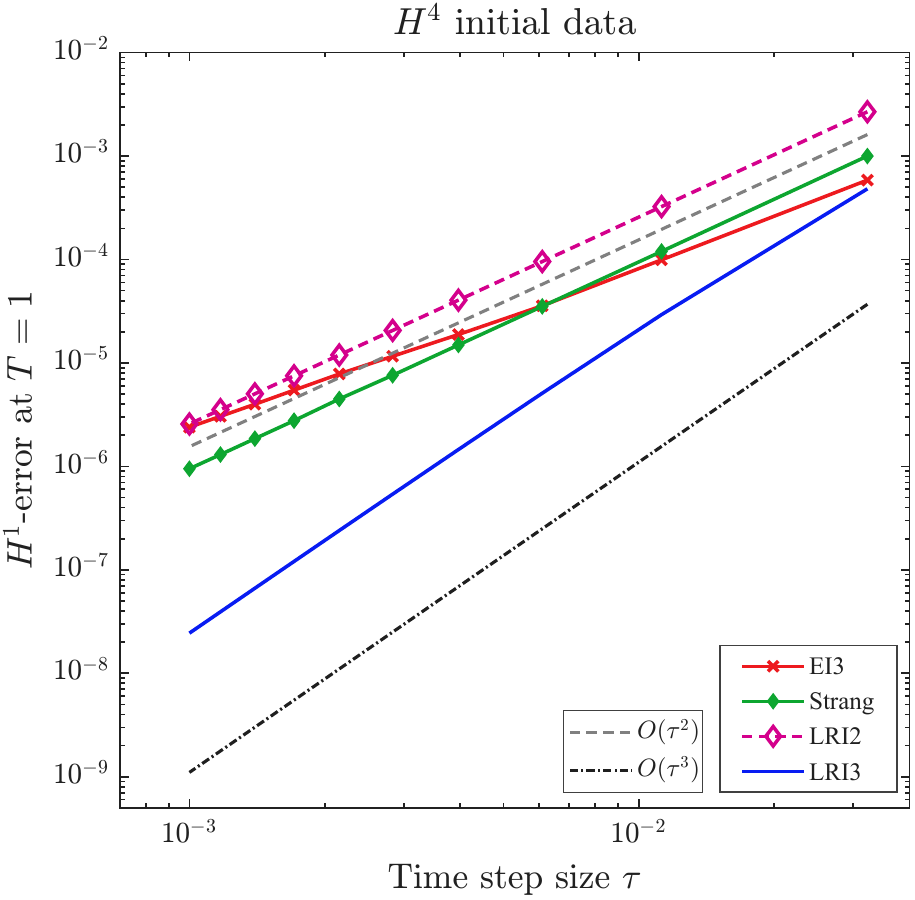}
\includegraphics[width=0.33\textwidth]{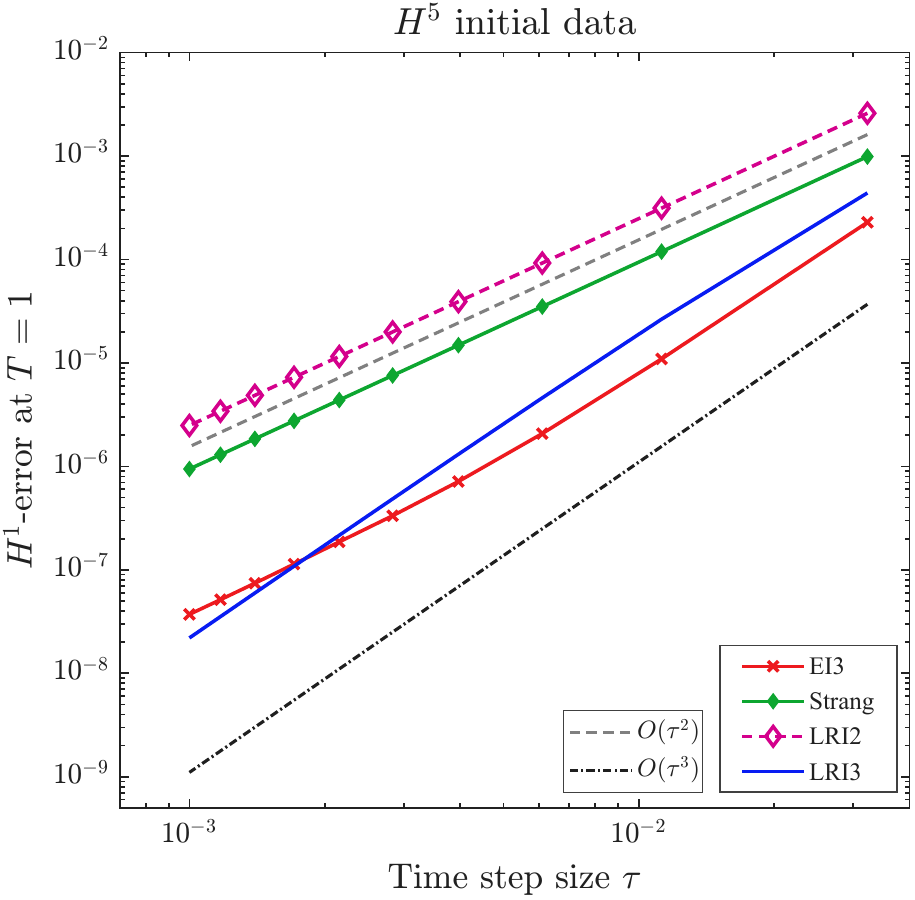}
\end{array}$$
\caption{$H^1$ errors of several time-stepping methods for initial data $u^0\in H^\theta$ with $\theta=3, 4, 5$.}
\label{fig:one}
\end{figure}

To exemplify the theoretical analysis in $H^r$ ($r>\tfrac12$),
we examine the errors measured in the $H^1$ norm.
Figure~\ref{fig:one} depicts these errors at $T=1$ for initial
data $u^0\in H^\theta$, with $\theta=3,4,$ and $5$.
The proposed method achieves third-order convergence for
$u^0\in H^3$, in agreement with \eqref{T1.1}, and maintains
this convergence rate for higher regularity.
By contrast, \texttt{EI3} still exhibits order reduction
when $u^0\in H^5$.
These results demonstrate the advantage of the proposed method
in achieving third-order accuracy under lower regularity.
\begin{figure}[htbp]
$$\begin{array}{cc}
\includegraphics[width=0.33\textwidth]{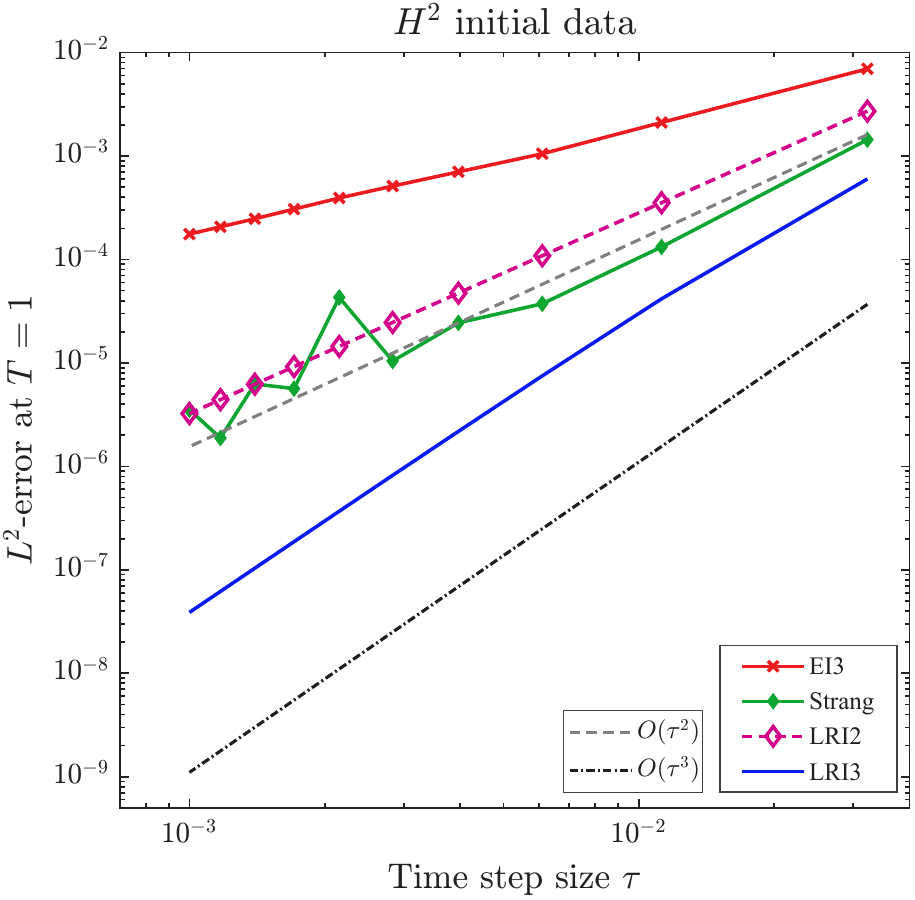}
\includegraphics[width=0.33\textwidth]{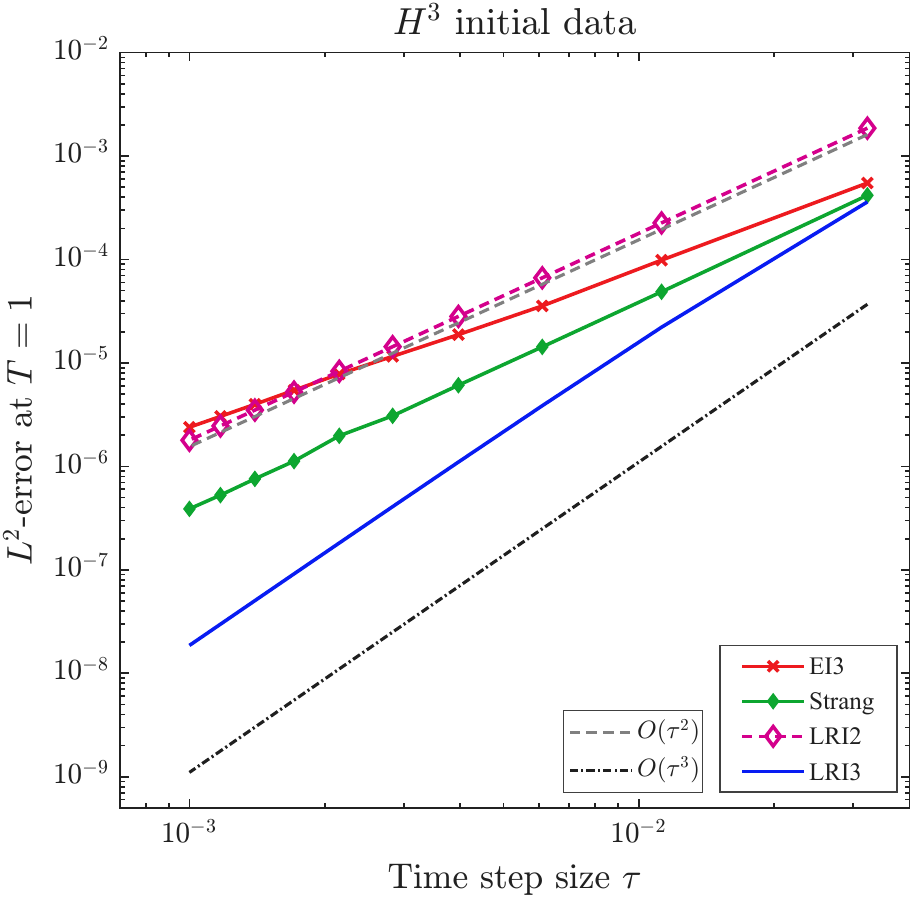}
\includegraphics[width=0.33\textwidth]{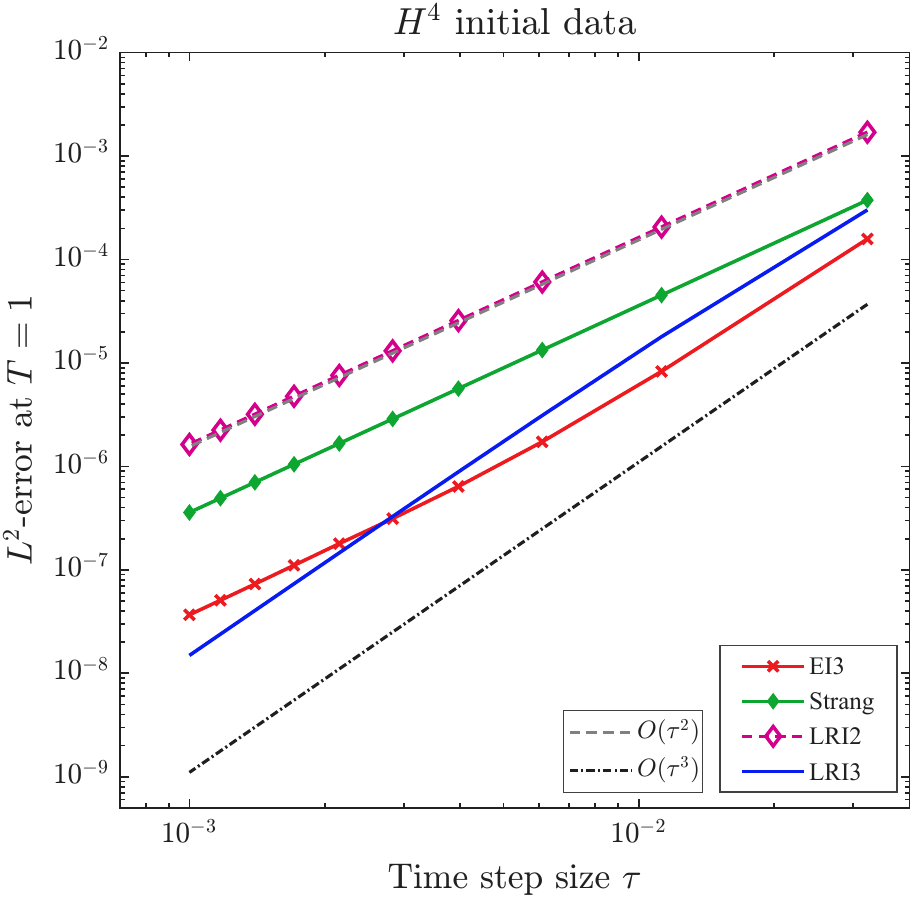}
\end{array}$$
\caption{$L^2$ errors of several time-stepping methods for initial data $u^0\in H^\theta$ with $\theta=2, 3, 4$.}
\label{fig:two}
\end{figure}
\par We further investigate the convergence behavior in $L^2$
for initial data $u^0\in H^\theta$ with $\theta=2,3,$ and $4$.
As illustrated in Figure~\ref{fig:two}, the convergence rate
of \texttt{LRI3} approaches three as the time step decreases,
even at the low regularity of $H^{2}$, corroborating
the theoretical estimate in \eqref{T1.2}.
Third-order convergence is consistently observed for smoother
initial data, while \texttt{EI3} still exhibits order reduction
up to $H^4$ regularity.
The numerical results thus confirm the improved low-regularity
convergence properties of \texttt{LRI3}.


In summary, the numerical results demonstrate that the proposed
\texttt{LRI3} scheme achieves third-order convergence in $H^1$
for initial data in $H^3$. Notably, third-order convergence in
$L^2$ is observed with initial regularity as low as $H^{2}$.
These findings confirm the theoretical results established
in Theorem~\ref{thm1.1}.

\section{Construction of  LRI3}\label{construction}
We now present the  formulation of \texttt{LRI3} for~\eqref{NLS}. For simplicity, we consider only the defocusing case $\lambda=1$, since the focusing case $\lambda=-1$ can be treated analogously. Following earlier works \cite{ostermann2018low,li2021fully,Wu2021,ostermann2022second,cao2024new,alama2024symmetric}, we introduce the twisted variable $v(t)=S_{-t}u(t)$
and apply Duhamel's formula to obtain
\begin{equation}\label{eq:3.1}
v(t_n+\tau)
= v(t_n)
- i \!\int_{0}^{\tau} S_{-(t_n+\sigma)}
\Big[
\big\lvert S_{t_n+\sigma}v(t_n+\sigma)\big\rvert^{2}\,
S_{t_n+\sigma}v(t_n+\sigma)
\Big]\mathrm{d}\sigma .
\end{equation}
The construction of the \texttt{LRI3} consists of the following three steps.\\
\noindent\textbf{Step 1.}
We use a second-order low-regularity method to approximate $v(t_n+\sigma)$ in \eqref{eq:3.1}, and after rearrangement obtain five integral terms to be treated.\\
\noindent\textbf{Step 2.}
For each of these five terms, we derive a numerical approximation
by exploiting the oscillatory structure of the phase, Taylor expansion,
and Fourier filtering. We then establish single-step and cumulative
remainder estimates. For the contributions requiring cancellation
across time steps, we employ discrete oscillatory estimates,
resonance counting, and discrete summation by parts.\\
\noindent\textbf{Step 3.}
By combining the numerical approximations of these five terms, we obtain the full scheme for the twisted variable $v(t)$. The final scheme for $u(t)$ then follows from the relation $u(t)=S_t v(t)$.
\par
More precisely, we first approximate $v(t_n+\sigma)$ by
\begin{equation}\label{ex}
v(t_n+\sigma)
= v(t_n)+\Omega^n(v(t_n),\sigma)+\mathcal{R}_{0,0}^{\sigma}(t_n),
\end{equation}
where  $
v(t_n)+\Omega^n(v(t_n),\sigma)=v(t_n)
+ \Omega_{1}^{n}(v(t_n),\sigma)
+ \Omega_{2}^{n}(v(t_n),\sigma)$ 
is the second-order low-regularity integrator constructed in~\cite{ostermann2022second}. Specifically,
\begin{subequations}\label{eq:3.4}
\begin{align}
\Omega_1^n(v(t_n),\sigma)&=-\ii\sigma S_{-t_n}\Big\{\big(S_{t_n}v(t_n)\big)^2\big(\varphi_1(-2\ii\sigma\partial_x^2)S_{-t_n}\overline{v}(t_n)\big)\Big\}\notag \\
&\quad\,\,+\ii\sigma S_{-(t_n+\sigma)}\Big\{\big(S_{t_n+\sigma}v(t_n)\big)^2\big(\psi_1(-2\ii\sigma\partial_x^2)S_{-(t_n-\sigma)}\overline{v}(t_n)\big)\Big\}\label{Omega1}\\ 
&\quad\,\,-\ii\sigma S_{-t_n}\Big\{\big(S_{t_n}v(t_n)\big)^2\big(\psi_1(-2\ii\sigma\partial_x^2)S_{-t_n}\overline{v}(t_n)\big)\Big\} \notag \\
&=:\Omega_{1,1}+\Omega_{1,2}+\Omega_{1,3}, \notag\\
\Omega_2^n(v(t_n),\sigma)&=-\frac{\sigma^2}{2}S_{-t_n}\Bigl\{|S_{t_n}v(t_n)|^4 S_{t_n}v(t_n)\Bigr\}. \label{Omega2}
\end{align}
\end{subequations}
Moreover, the following estimate for the remainder term $\mathcal{R}_{0,0}^{\sigma}(t_n)$ was proved in \cite{ostermann2022second}:
\begin{equation}\label{eq:R00}
\bigl\|\mathcal{R}_{0,0}^{\sigma}(t_n)\bigr\|_{H^{r-2}}
\le  C\!\left(\,\|v\|_{C([0,T];H^{r})}\right)\sigma^{3}\quad \text{for any } r>\frac{5}{2},
\end{equation}
where $C\!\left(\,\|v\|_{C([0,T];H^{r})}\right)$ denotes a constant depending only on
$\|v\|_{C([0,T];H^{r})}$.
\par Inserting~\eqref{ex} into the right-hand of~\eqref{eq:3.1}, we derive that 
\begin{equation}\label{integral}
\begin{aligned}
v(t_n+\tau)
=&v(t_n)
-\ii\int_{0}^{\tau}\!\!
S_{-(t_n+\sigma)}
\Big[
S_{t_n+\sigma}
\Big(v(t_n)+\sum\limits_{j=1}^2\Omega_{j}^{n}(v(t_n),\sigma)+\mathcal{R}_{0,0}^{\sigma}(t_n)\Big)
\\[-0.2em]
&\qquad\qquad \qquad \,\,\, \,\cdot
\Big\lvert
S_{t_n+\sigma}
\Big(v(t_n)+\sum\limits_{j=1}^2\Omega_{j}^{n}(v(t_n),\sigma)+\mathcal{R}_{0,0}^{\sigma}(t_n)\Big)
\Big\rvert^{2}
\Big]\,\mathrm{d}\sigma .
\end{aligned}
\end{equation} 
By a straightforward computation, we obtain the following expression:
\begin{equation}\label{expressions}
	v(t_n+\tau)
	= v(t_n)
	+ I_{(3)}^\tau(t_n)
	+ I_{(5)}^\tau(t_n)
	+ I_{(7)}^\tau(t_n)
	+ \mathcal{R}_{(0)}^\tau(t_n).
\end{equation}
with the notation
\begin{align*} 	
I_{(3)}^\tau(t_n):=& -\ii\int_{0}^{\tau} S_{-(t_n+\sigma)}\Big[ \big|S_{t_n+\sigma}v(t_n)\big|^2 \big(S_{t_n+\sigma}v(t_n)\big) \Big] \mathrm{d}\sigma,\\
I_{(5)}^\tau(t_n):=&-2\ii\int_{0}^{\tau} 
S_{-(t_n+\sigma)}\Big[ \big|S_{t_n+\sigma}v(t_n)\big|^2 \big(S_{t_n+\sigma}\,\Omega_1^n(v(t_n),\sigma)\big) \Big] \mathrm{d}\sigma\\
&-\ii\int_{0}^{\tau}
S_{-(t_n+\sigma)}\Big[ \big(S_{t_n+\sigma}v(t_n)\big)^2 \big(S_{-(t_n+\sigma)}\,\overline{\Omega_1^n(v(t_n),\sigma)}\big) \Big] \mathrm{d}\sigma,\\
I_{(7)}^\tau(t_n):=&-2\ii\int_{0}^{\tau}  S_{-(t_n+\sigma)}\Big[ \big|S_{t_n+\sigma}v(t_n)\big|^2 \big(S_{t_n+\sigma}\,\Omega_2^n(v(t_n),\sigma)\big) \Big] \mathrm{d}\sigma\\
&-\ii\int_{0}^{\tau} S_{-(t_n+\sigma)}\Big[ \big(S_{t_n+\sigma}v(t_n)\big)^2 \big(S_{-(t_n+\sigma)}\,\overline{\Omega_2^n(v(t_n),\sigma)}\big) \Big] \mathrm{d}\sigma\\
&-2\ii\int_{0}^{\tau}  S_{-(t_n+\sigma)}\Big[ \big|S_{t_n+\sigma}\,\Omega_1^n(v(t_n),\sigma)\big|^2 \big(S_{t_n+\sigma}v(t_n)\big) \Big] \mathrm{d}\sigma\\
&-\ii\int_{0}^{\tau} S_{-(t_n+\sigma)}\Big[\big(S_{t_n+\sigma}\,\Omega_1^n(v(t_n),\sigma)\big)^2 \big(S_{-(t_n+\sigma)}\,\overline{v}(t_n)\big) \Big] \mathrm{d}\sigma.
\end{align*}
Moreover, the remainder $\mathcal{R}_{(0)}^\tau(t_n)$ satisfies the estimates below,
and the proof of these estimates will be deferred to Subsection~\ref{pooflemma3.1}.
\begin{lemma}\label{lemma3.1}
Let $0<\tau\le1$. Then the following bounds hold.
	\begin{enumerate}
		\item If $r>\tfrac52$ and $v\in C([0,T];H^r)$, then
		\[
		\|\mathcal{R}_{(0)}^\tau(t_n)\|_{H^{r-2}}
		\le C\big(\|v\|_{C([0,T];H^r)}\big)\tau^4,\qquad 0\le n < L.
		\]
		
		\item If $v\in C([0,T];H^2)$, then
		\[
		\max_{1\le n\le L}
		\left\|\sum_{m=0}^{n-1}\mathcal{R}_{(0)}^\tau(t_m)\right\|_{L^2}
		\le C_T\tau^3,
		\]
		where $C_T$ depends only on $T$ and $\|v\|_{C([0,T];H^2)}$.
	\end{enumerate}
\end{lemma}

In the remainder of this section, we derive the numerical approximations for all integrals appearing in \eqref{expressions}. Based on these approximations, the \texttt{LRI3} scheme is formulated.

\subsection{Approximation of the cubic term $I_{(3)}^\tau(t_n)$}
For the cubic contribution $I_{(3)}^\tau(t_n)$, we use the
approximation of \cite{cao2024new} and write
\[
I_{(3)}^\tau(t_n)
= \Psi_{(3)}^n(v(t_n),\tau)
+ \mathcal{R}_{(3)}^\tau(t_n),
\]
where $\Psi_{(3)}^n(v(t_n),\tau)$ is their cubic approximation
$\Psi_1^n(v(t_n),\tau)$.
To express the remainder $\mathcal{R}_{(3)}^\tau(t_n)$ in Fourier
space, we define
$
\mathcal{M}_\tau(\alpha)
= \frac{1}{\tau}\int_0^\tau
\sigma\bigl(\mathrm e^{2\ii\sigma\alpha}-1\bigr)\,\mathrm{d}\sigma.
$
Throughout the remainder of the paper, all Fourier indices take
integer values. For brevity, we write
$
\{\hat{\overline{v}}_{k_1}\hat{v}_{k_2}\hat{v}_{k_3}\}(t_n)
:= \hat{\overline{v}}_{k_1}(t_n)
\hat{v}_{k_2}(t_n)\hat{v}_{k_3}(t_n)
$.
With this notation, the cubic remainder is given by
\begin{equation}\label{R1tau}
	\mathcal R_{(3)}^\tau(t_n)
	=\ii\sum_{k\ne0}\sum_{k_1+k_2+k_3=k}
	\mathrm e^{\ii t_n\phi}
	\bigl\{\hat{\overline v}_{k_1}\hat v_{k_2}\hat v_{k_3}\bigr\}(t_n)
	\mathrm e^{\ii kx} \cdot K_{3,R}^\tau,
\end{equation}
where $A=kk_1$, $B=k_2k_3$, and
$\phi=k^2+k_1^2-k_2^2-k_3^2=2(A+B)$.
Writing $D_X(\sigma)=e^{2\ii\sigma X}-1$, the cubic remainder
kernel can be written as
\begin{equation*}
	K_{3,R}^\tau
	=\frac{k_1}{k}\int_0^\tau
	D_A(\sigma)
	\left[D_B(\sigma)-\frac{\sigma}{\tau}D_B(\tau)\right]
	\,\mathrm d\sigma
   \quad \text{for }\,k\ne0.
\end{equation*}
The zero Fourier mode of the remainder vanishes.
Although the cubic approximation is taken from \cite{cao2024new},
our error analysis exploits discrete oscillatory cancellation across
time steps to establish a cumulative remainder estimate under
weaker regularity assumptions.
The following lemma provides both single-step and cumulative
estimates for $\mathcal R_{(3)}^\tau(t_n)$.
The proof is given in Subsection~\ref{pooflemma3.2}.
\begin{lemma}\label{lemma3.2}
	The following estimates hold for $0< \tau \le1$:
	\begin{enumerate}
		\item If $r>\tfrac52$ and $v\in C([0,T];H^r)$, then
		\[
		\|\mathcal R_{(3)}^\tau(t_n)\|_{H^{r-2}}
		\lesssim
		\tau^4\|v\|_{C([0,T];H^r)}^3,\qquad 0\le n < L.
		\]
		
		\item If $v\in C([0,T];H^{2})$, then
		\[
		\max_{0\le n\le L}
		\left\|\sum_{m=0}^{n-1}
		\mathcal R_{(3)}^\tau(t_m)\right\|_{L^2}
		\le C\tau^3,
		\]
		where $C$ depends only on $T$ and
		$\|v\|_{C([0,T];H^{2})}$.
	\end{enumerate}
\end{lemma}

\subsection{Approximation of the quintic term $I_{(5)}^\tau(t_n)$}
\label{subsec:quintic-approximation}
Let $\boldsymbol{k}=(k_1,\ldots,k_5)$. Substituting \eqref{Omega1} into the two
integrals defining $I_{(5)}^\tau(t_n)$ in \eqref{expressions}
gives
\begin{equation}\label{eq:quintic-exact}
	I_{(5)}^\tau(t_n)
	=\mathcal T_{1,n}[-2K_{5,1}^{1,\tau}+2K_{5,1}^{2,\tau}]\big(v(t_n)\big)
	+\mathcal T_{2,n}[K_{5,2}^{1,\tau}-K_{5,2}^{2,\tau}]\big(v(t_n)\big).
\end{equation}
Here $p=1,2$ refers to the first and second integrals, respectively.
For a kernel $K=K(\boldsymbol{k})$, define
\[
\widehat{\mathcal T_{p,n}[K]\big(w(t_n)\big)}_k
=\sum_{k_1+\cdots+k_5=k}
\mathrm e^{\ii t_n\phi_p}K(\boldsymbol{k})\,V_p\big(w(t_n);\boldsymbol{k}\big),
\]
where
\begin{align*}
	V_1\big(w(t_n);\boldsymbol{k}\big)
	&=\{\widehat{\overline w}_{k_1}\widehat w_{k_2}\widehat w_{k_3}
	\widehat w_{k_4}\widehat{\overline w}_{k_5}\}(t_n),
	&\phi_1&=k^2+k_1^2-k_2^2-k_3^2-k_4^2+k_5^2,\\
	V_2\big(w(t_n);\boldsymbol{k}\big)
	&=\{\widehat w_{k_1}\widehat w_{k_2}\widehat{\overline w}_{k_3}
	\widehat{\overline w}_{k_4}\widehat w_{k_5}\}(t_n),
	&\phi_2&=k^2-k_1^2-k_2^2+k_3^2+k_4^2-k_5^2.
\end{align*}

For $p=1,2$, set $a_p=\alpha_p^2$ and $b=k_5^2$, where
$\alpha_1=k_1$ and $\alpha_2=k_3+k_4$. Define
\begin{align*}
	\beta_1 &= k_1k_2+(k_1+k_2)(k_3+k_4+k_5),
	& \beta_2 &= \beta_1+2(k_3+k_4)k_5,\\
	\beta_3 &= k_3k_4+(k_3+k_4)k_5,
	& \Gamma_1 &= \beta_1+\beta_3,\qquad \Gamma_2=\beta_2-\beta_3.
\end{align*}
To express the kernels in a common form, we introduce
\[
F(a,b) = \int_0^\tau s \mathrm e^{2\ii sa}\varphi_1(2\ii sb)\,\mathrm ds,
\qquad
Q_p(a,b) = \int_0^\tau s^2 \mathrm e^{2\ii sa}q_p(2\ii sb)\,\mathrm ds,
\]
where $q_1=\psi_1$ and $q_2=\varphi_2$. We also set
$G(a,b)=\partial_a F(a,b)$.
Using $\varphi_1(-z)=\mathrm e^{-z}\varphi_1(z)$ and
$\psi_1(-z)=-\mathrm e^{-z}\varphi_2(z)$, we write the four exact kernels
in \eqref{eq:quintic-exact} as
\begin{align}
	K_{5,p}^{1,\tau}
	= F(a_p+\beta_p,b), \qquad K_{5,p}^{2,\tau}
	= \int_0^\tau s \mathrm e^{2\ii sa_p}
	\Theta_p(s)q_p(2\ii sb)\,\mathrm ds,
\end{align}
where
\[
\Theta_1(s)=\mathrm e^{2\ii s\Gamma_1}(1-\mathrm e^{-2\ii s\beta_3}),
\qquad
\Theta_2(s)=-\mathrm e^{2\ii s\Gamma_2}(1-\mathrm e^{2\ii s\beta_3}).
\]
Expanding $K_{5,p}^{1,\tau}$ to first order in $\beta_p$ gives
$F(a_p,b)+\beta_p G(a_p,b)$, since $G=\partial_a F$.
For $K_{5,p}^{2,\tau}$, replacing $\Theta_p(s)$ by its linear term
$2\ii s\beta_3$ gives $2\ii\beta_3 Q_p(a_p,b)$.
We then approximate each of $F$, $G$, and $Q_p$ by
\[
M^{\mathrm{approx}}(a,b)=M(a,0)+M(0,b)-M(0,0),
\qquad M\in\{F,G,Q_p\}.
\]
Applying the filter
$\mathcal P_\tau(\boldsymbol k)=\prod_{\ell=1}^5(1+\tau|k_\ell|)^{-1}$
to the linear corrections gives
\begin{equation}\label{eq:quintic-approximate-kernels}
	\begin{aligned}
		\widetilde K_{5,p}^{1,\tau}
		=F^{\mathrm{approx}}(a_p,b)
		+\beta_p\mathcal P_\tau G^{\mathrm{approx}}(a_p,b),\quad\widetilde K_{5,p}^{2,\tau}
		=2\ii\beta_3\mathcal P_\tau Q_p^{\mathrm{approx}}(a_p,b),
		\quad p=1,2.
	\end{aligned}
\end{equation}
Let $P_\tau$ denote the Fourier multiplier defined by
$\widehat{P_\tau f}_k=(1+\tau|k|)^{-1}\hat f_k$.
Writing $u_\tau=P_\tau u(t_n)$ and $\overline{u_\tau}=P_\tau\bar{u}(t_n)$, with
$u(t_n)=S_{t_n}v(t_n)$, we substitute the approximate kernels into
\eqref{eq:quintic-exact} to obtain
\begin{align*}
	\Psi_{(5)}^n(v(t_n),\tau)
	&=2\tau^2 S_{-t_n}\Big\{
	u(t_n)^2|u(t_n)|^2
	\Big[\big(\psi_{1,2}^-+\tfrac12-\varphi_{2,2}^-\big)
	\overline{u}(t_n)\Big]\Big\}-\tau^2 S_{-t_n}\Big\{
	u(t_n)^3\big(\psi_{1,2}^-+\tfrac12\big)
	\big(\overline{u}(t_n)^2\big)\Big\}\\
	&\quad+\tau^2 S_{-t_n}\Big\{
	|u(t_n)|^4\big(\varphi_{2,2}^-u(t_n)\big)\Big\}+4\ii\tau^3 S_{-t_n}\Big\{
	|u_\tau|^2(Du_\tau)^2
	\Big[\big(\psi_{2,2}^--\tfrac16+\Phi_{2,2}^-\big)\overline{u_\tau}\Big]\Big\}\\
	&\quad+4\ii\tau^3 S_{-t_n}\Big\{
	D(|u_\tau|^2)\,D\big[u_\tau^2(\Phi_{1,2}^-\overline{u_\tau})\big]\Big\}+8\ii\tau^3 S_{-t_n}\Big\{
	u_\tau|u_\tau|^2(Du_\tau)
	\Big[\big(\psi_{2,2}^--\tfrac16+\Phi_{2,2}^-\big)D\overline{u_\tau}\Big]\Big\}\\
	&\quad+8\ii\tau^3 S_{-t_n}\Big\{
	D\Big[u_\tau\Big(\big(\psi_{2,2}^--\tfrac16\big)\overline{u_\tau}\Big)\Big]\,
	D(u_\tau|u_\tau|^2)\Big\}+4\ii\tau^3 S_{-t_n}\Big\{
	u_\tau^2(Du_\tau)(D\overline{u_\tau})(\Phi_{1,2}^-\overline{u_\tau})\Big\}\\
	&\quad+8\ii\tau^3 S_{-t_n}\Big\{
	u_\tau^2(Du_\tau)(D\overline{u_\tau})\Big[\big(\psi_{2,2}^--\tfrac16\big)\overline{u_\tau}\Big]\Big\}-12\ii\tau^3 S_{-t_n}\Big\{
	u_\tau^2(Du_\tau)\big(\psi_{2,2}^--\tfrac16\big)
	\big[\overline{u_\tau}(D\overline{u_\tau})\big]\Big\}\\
	&\quad+4\ii\tau^3 S_{-t_n}\Big\{
	u_\tau|u_\tau|^2(D\overline{u_\tau})
	\Big[\big(\varphi_{3,2}^--2\Phi_{1,2}^-\big)Du_\tau\Big]\Big\}\\
	&\quad-2\ii\tau^3 S_{-t_n}\Big\{
	D(u_\tau^2)\,D\Big[
	2u_\tau\big(\psi_{2,2}^--\tfrac16\big)(\overline{u_\tau}^{\,2})
	+\overline{u_\tau}^{\,2}(\Phi_{1,2}^-u_\tau)
	\Big]\Big\}\\
	&\quad-2\ii\tau^3 S_{-t_n}\Big\{
	(Du_\tau)^2\Big[
	2u_\tau\big(\psi_{2,2}^--\tfrac16\big)(\overline{u_\tau}^{\,2})
	+\overline{u_\tau}^{\,2}(\Phi_{1,2}^-u_\tau)
	\Big]\Big\}\\
	&\quad+2\ii\tau^3 S_{-t_n}\Big\{
	u_\tau^3\big(\psi_{2,2}^--\tfrac16\big)\big[(D\overline{u_\tau})^2\big]
	+u_\tau^2(D\overline{u_\tau})^2(\varphi_{3,2}^-u_\tau)\Big\}.
\end{align*}
The exact quintic term can thus be written as
\begin{equation}\label{eq:quintic-decomposition}
	I_{(5)}^\tau(t_n)
	= \Psi_{(5)}^n(v(t_n),\tau) + \mathcal R_{(5)}^\tau(t_n).
\end{equation}
The Fourier representation of $\mathcal R_{(5)}^\tau(t_n)$ is
obtained from \eqref{eq:quintic-exact} by replacing each kernel
$K_{5,p}^{j,\tau}$ with the remainder kernel
$K_{5,p,R}^{j,\tau}=K_{5,p}^{j,\tau}-\widetilde K_{5,p}^{j,\tau}$.
The following lemma gives both single-step and cumulative bounds for
the quintic remainder. Its proof is deferred to
Subsection~\ref{sec:quintic-remainder-proof}.
\begin{lemma}\label{lem:quintic-remainder}
	The following estimates hold for $0<\tau\le 1$:
	\begin{enumerate}
		\item If $r>\tfrac52$ and $v\in C([0,T];H^r)$, then
		\[
		\|\mathcal R_{(5)}^\tau(t_n)\|_{H^{r-2}}
		\lesssim
		\tau^4\|v\|_{C([0,T];H^r)}^5,\qquad 0\le n < L.
		\]
		\item If
		$v\in C([0,T];H^{2})$, then
		\[
		\max_{0\le n\le L}
		\left\|\sum_{m=0}^{n-1}
		\mathcal R_{(5)}^\tau(t_m)\right\|_{L^2}
		\le C\tau^3,
		\]
		where $C$ depends on $T$ and $\|v\|_{C([0,T];H^2)}$,
		and is independent of $\tau$.
	\end{enumerate}
\end{lemma}
\subsection{Approximation of the septic term $I_{(7)}^\tau(t_n)$}
Let $I_{(7)}^{j,\tau}(t_n)$ denote the $j$th summand in the
definition of $I_{(7)}^\tau(t_n)$, for $j=1,\ldots,4$.
In the Fourier representation of each summand, we retain the
principal oscillation and evaluate the remaining oscillatory
factors at $\sigma=0$.
To illustrate this procedure, consider the first integral. Its phase splits as
\[
k^2+k_1^2-k_2^2-(k_3+\cdots+k_7)^2
=2k_1^2+2\beta_4,
\qquad
\beta_4=k_1k_2+(k_1+k_2)(k_3+\cdots+k_7).
\]
Replacing $\mathrm e^{2\ii\sigma\beta_4}$ by $1$ gives
$
\int_0^\tau \sigma^2 \mathrm e^{2\ii\sigma k_1^2}\,\mathrm d\sigma
=2\tau^3\psi_2(2\ii\tau k_1^2).
$
Applying the same procedure to the remaining integrals gives
the following contributions. The last column lists the
physical-space approximations after factoring out
$\ii\tau^3S_{-t_n}$:
\[
\begingroup
\renewcommand{\arraystretch}{1.3}
\begin{array}{c|c|l}
	j & \text{Retained phase} & \text{Contribution}\\ \hline
	1 & 2k_1^2
	& 2(S_{t_n}v(t_n))^2 |S_{t_n}v(t_n)|^4
	\bigl(\psi_{2,2}^- S_{-t_n}\overline{v}(t_n)\bigr)\\
	2 & 2(k_3+\cdots+k_7)^2
	& (S_{t_n}v(t_n))^2
	\psi_{2,2}^-\big\{|S_{t_n}v(t_n)|^4
	S_{-t_n}\overline{v}(t_n)\big\}\\
	3 & 2(k_1+k_2+k_3)^2
	& -4(S_{t_n}v(t_n))^2 |S_{t_n}v(t_n)|^2
	\psi_{2,2}^-\big\{|S_{t_n}v(t_n)|^2
	S_{-t_n}\overline{v}(t_n)\big\}\\
	4 & 2k_7^2
	& 2(S_{t_n}v(t_n))^2 |S_{t_n}v(t_n)|^4
	\bigl(\psi_{2,2}^- S_{-t_n}\overline{v}(t_n)\bigr)
\end{array}
\endgroup
\]
Combining the first and fourth contributions, we obtain
\begin{align*}
	\Psi_{(7)}^n(v(t_n),\tau)
	&= 4\ii\tau^3 S_{-t_n}\Big\{(S_{t_n}v(t_n))^2
	\big[
	|S_{t_n}v(t_n)|^4
	\bigl(\psi_{2,2}^- S_{-t_n}\overline{v}(t_n)\bigr)
	+\frac 14\psi_{2,2}^-\bigl(|S_{t_n}v(t_n)|^4
	S_{-t_n}\overline{v}(t_n)\bigr)\\
	&\qquad\qquad\quad
	-|S_{t_n}v(t_n)|^2
	\psi_{2,2}^-\big\{|S_{t_n}v(t_n)|^2
	S_{-t_n}\overline{v}(t_n)\big\}
	\big]\Big\}.
\end{align*}
Thus, the septic term decomposes as
\[
I_{(7)}^\tau(t_n)
= \Psi_{(7)}^n(v(t_n),\tau) + \mathcal R_{(7)}^\tau(t_n).
\]
The following lemma bounds the remainder $\mathcal R_{(7)}^\tau(t_n)$.
Its proof is given in Subsection~\ref{pooflemma3.4}.
\begin{lemma}\label{lemma3.4}
	Let $r\ge2$ and $v\in C([0,T];H^r)$. Then we have the following estimate:
	\begin{equation*}
		\begin{aligned}
			\| \mathcal{R}_{(7)}^\tau(t_n)\|_{H^{r-2}}\lesssim\tau^4\| v\|_{C([0,T];H^r)}^7,\qquad 0\le n<L.
		\end{aligned}
	\end{equation*}
\end{lemma}

Overall, by combining all the above cases, we obtain 
\begin{equation}\label{expression}
	v(t_{n+1})=v(t_n)+\sum_{j\in\{3,5,7\}}\Psi_{(j)}^n(v(t_n),\tau)+\sum_{j\in\{0,3,5,7\}}\mathcal R_{(j)}^\tau(t_n)=:\Psi^n(v(t_n),\tau)
	+\mathcal R^\tau(t_n).
\end{equation}
The numerical scheme is defined by omitting $\mathcal{R}^\tau(t_n)$, i.e.,
\begin{equation}\label{vscheme}
v^{n+1}=\Psi^n(v^n,\tau),\quad n\ge0;\quad \text{for}\,\,\,\,v^0=u^0.
\end{equation}
Substituting $v^n = S_{-t_n}u^n$ into \eqref{vscheme} yields 
$u^{n+1}=S_{t_{n+1}}\Psi^n(S_{-t_n}u^n,\tau),$
which gives  the \texttt{LRI3} scheme \eqref{littlescheme}-\eqref{numericalscheme}.

\section{Convergence analysis (proof of Theorem~\ref{thm1.1})}\label{proof}
In Subsection~\ref{p1}, we prove the technical lemmas stated
in Section~\ref{construction}. In Subsection~\ref{r12},
we establish the local error and stability estimates in
$H^r$ for $r>\tfrac12$ to prove \eqref{T1.1}.
In Subsection~\ref{r13}, we first establish the cumulative
local error bound in $L^2$ and derive an a priori $H^1$ bound
for the numerical solution, which is then used to obtain an
$L^2$ stability estimate. Combining the cumulative local error
bound with the $L^2$ stability estimate yields \eqref{T1.2},
thereby completing the proof of Theorem~\ref{thm1.1}.
\subsection{Proofs of Lemmas \ref{lemma3.1}-\ref{lemma3.4}}\label{p1}
This subsection proves Lemmas~\ref{lemma3.1}--\ref{lemma3.4} stated in Section~\ref{construction}. The proof of Lemma~\ref{lemma3.1} is deferred to the end
of this subsection, as it uses an argument similar to
that in the proof of Lemma~\ref{lem:quintic-remainder}.
\subsubsection{Proof of Lemma~\ref{lemma3.2}}\label{pooflemma3.2}
\begin{proof}
	Recall that $A=kk_1$, $B=k_2k_3$, and $\phi=2(A+B)$.
	For $0\le \sigma\le\tau$, the factors in $K_{3,R}^\tau$ satisfy
	$
	|D_A(\sigma)|\lesssim \min\{1,\tau|A|\},
	\,
	\left|D_B(\sigma)-\frac{\sigma}{\tau}D_B(\tau)\right|
	\lesssim \min\{1,(\tau B)^2\}.
	$
	Since
	$|k_1A/k|=|k_1|^2$ for $k\ne0$ and the remainder vanishes at
	$k=0$, we obtain
	\begin{equation}\label{eq:cubic-remainder-kernel-bound}
		|K_{3,R}^\tau|
		\lesssim
		\frac{\tau^4|k_1k_2k_3|^2}
		{(1+\tau|A|)(1+(\tau B)^2)}.
	\end{equation}
	Therefore, by Lemma~\ref{kp}~\textup{(i)}, we obtain, for any
	$r>\tfrac{5}{2}$,
	\[
	\|\mathcal R_{(3)}^\tau(t_n)\|_{H^{r-2}}
	\lesssim \tau^4 \|v\|_{C([0,T];H^r)}^3.
	\]
	
	For the cumulative estimate, let $\mathcal D_m(f_1,f_2,f_3)$
	be the operator obtained from \eqref{R1tau} by setting $n=m$
	and replacing $\overline v(t_m),v(t_m),v(t_m)$ with
	$\overline{f_1},f_2,f_3$, respectively,
	and set $\mathcal A_n=\sum_{m=0}^{n-1}\mathcal D_m$.
	For fixed inputs, Lemma~\ref{lem:discrete-time-summation}(i)
	and \eqref{eq:cubic-remainder-kernel-bound} show that, after
	factoring out $\tau^3$ and the three input weights
	$\langle k_j\rangle^2$, the Fourier multiplier
	of $\mathcal A_n$ is bounded in modulus by $C_T W_\tau$,
	uniformly in $n$, where
	\[
	W_\tau
	=\frac{g_\tau(\phi)}
	{(1+\tau|A|)(1+(\tau B)^2)}.
	\]
By $L^2$ duality, to estimate
$\|\mathcal A_n(f_1,f_2,f_3)\|_{L^2}$,
it suffices to bound the following sum:
	\[
	\sum_{k\in\mathbb Z}\sum_{k_1+k_2+k_3=k}
	W_\tau\,h_0(k)h_1(k_1)h_2(k_2)h_3(k_3),
	\]
	where $h_0=|\hat q|$ for a unit $L^2$ test function $q$,
	and $h_1,h_2,h_3$ are the absolute Fourier coefficients of
	$\overline{f_1},f_2,f_3$ weighted by $\langle k_j\rangle^2$.
	Applying Cauchy--Schwarz to
	$W_\tau^{3/8}h_0h_1$ and $W_\tau^{5/8}h_2h_3$,
	we bound their squared sums by
	$M_{01}\|h_0\|_{\ell^2}^2\|h_1\|_{\ell^2}^2$ and
	$M_{23}\|h_2\|_{\ell^2}^2\|h_3\|_{\ell^2}^2$, respectively, where
	\[
	M_{01}=\sup_{k,k_1}
	\sum_{k_2+k_3=k-k_1}W_\tau^{3/4},
	\quad
	M_{23}=\sup_{k_2,k_3}
	\sum_{k_1\in\mathbb Z}W_\tau^{5/4},
	\]
	with $k=k_1+k_2+k_3$ in the second sum. Restoring the factor
	$\tau^3$ from the multiplier bound gives
	\[
	\|\mathcal A_n(f_1,f_2,f_3)\|_{L^2}
	\lesssim_T \tau^3(M_{01}M_{23})^{1/2}
	\prod_{j=1}^3\|f_j\|_{H^2}.
	\]
	
	For $M_{01}$, $A$ is fixed and $\phi$ is quadratic in $k_2$
	with leading coefficient $-2$. On the resonance interval
	$I_\ell$ of Lemma~\ref{lem:discrete-time-summation}(ii),
	$\tau B=\pi\ell-\tau A+O(1)$. Applying that lemma with
	$\rho=\tfrac34$, we obtain
	\[
	\begin{aligned}
		M_{01}
		\lesssim
		\sup_{k,k_1}\sum_{\ell\in\mathbb Z}
		\frac{1}{(1+|\pi\ell-\tau A|^2)^{3/4}}
		\sum_{\substack{k_2+k_3=k-k_1\\ \phi\in I_\ell}}
		g_\tau(\phi)^{3/4}\lesssim
		\sup_A\sum_{\ell\in\mathbb Z}
		\frac{1}{(1+|\pi\ell-\tau A|^2)^{3/4}}
		\lesssim 1.
	\end{aligned}
	\]
	
	For $M_{23}$, $B$ is fixed and $\phi$ is quadratic in $k_1$
	with leading coefficient $2$. On each resonance interval
	$I_\ell$, we have $\tau A=\pi\ell-\tau B+O(1)$.
	Applying Lemma~\ref{lem:discrete-time-summation}(ii) with
	$\rho=\tfrac54$, we obtain
	\[
	\begin{aligned}
		M_{23}
		\lesssim
		\sup_{k_2,k_3}\sum_{\ell\in\mathbb Z}
		\frac{1}{(1+|\pi\ell-\tau B|)^{5/4}}
		\sum_{\substack{k_1\in\mathbb Z\\ \phi\in I_\ell}}
		g_\tau(\phi)^{5/4}\lesssim
		\sup_B\sum_{\ell\in\mathbb Z}
		\frac{1}{(1+|\pi\ell-\tau B|)^{5/4}}
		\lesssim 1.
	\end{aligned}
	\]
	Consequently,
	$
	\|\mathcal A_n(f_1,f_2,f_3)\|_{L^2}
	\lesssim_T \tau^3\prod_{j=1}^3\|f_j\|_{H^2}.
	$
By the twisted equation, we have
	\[
	\sum_{m=0}^{L-2}
	\|v(t_{m+1})-v(t_m)\|_{H^2}
	\lesssim T\|v\|_{C([0,T];H^2)}^3.
	\]
	The cumulative estimate follows from
	Lemma~\ref{lem:discrete-time-summation}(iii) with
	$d=3$, $X=H^2$, $Y=L^2$, and $w_m=v(t_m)$.
\end{proof}
\subsubsection{Proof of Lemma~\ref{lem:quintic-remainder}}\label{sec:quintic-remainder-proof}
\begin{proof}
	Set $\mu(x)=\min\{1,x\}$ for $x\ge0$. For brevity, we write $F$, $G$, and $Q_p$ for $F(a_p,b)$,
	$G(a_p,b)$, and $Q_p(a_p,b)$, respectively, and adopt the
	same convention for their approximations. Expanding to first order in $\beta_p$ and $\beta_3$ and applying
	Lemma~\ref{jb} to the resulting exponential differences and
	Taylor remainders, we obtain
	\begin{align*}
		|F(a_p+\beta_p,b)-F-\beta_pG|
		\lesssim \frac{\tau^4|\beta_p|^2}{1+\tau|\beta_p|},\quad
		|K_{5,p}^{2,\tau}-2\ii\beta_3Q_p|
		\lesssim \frac{\tau^4|\beta_3|^2}{1+\tau|\beta_3|}
		+\tau^2\mu(\tau|\beta_3|)\mu(\tau|\Gamma_p|).
	\end{align*}
	The separable approximation errors satisfy
	\begin{align*}
		|F-F^{\mathrm{approx}}|
		\lesssim \tau^2\mu(\tau a_p)\mu(\tau b),\quad|G-G^{\mathrm{approx}}|+|Q_p-Q_p^{\mathrm{approx}}|
		\lesssim \tau^3\mu(\tau a_p)\mu(\tau b),
	\end{align*}
	since their integrands contain
	$(\mathrm e^{2\ii sa_p}-1)(q(2\ii sb)-q(0))$ with
	$q\in\{\varphi_1,\psi_1,\varphi_2\}$. The filter errors are bounded by
	$C\tau^4(|\beta_p|+|\beta_3|)\sum_{\ell=1}^5|k_\ell|$.
	
	Let $\mathscr I$ consist of the multi-indices
	$\boldsymbol d\in\{0,1,2\}^5$ with $\sum_{\ell=1}^5 d_\ell\le4$
	and at most one entry equal to two. The product of two distinct
	frequency pairs $|k_i k_q|$ and $|k_u k_v|$ gives a monomial
	controlled by $\prod_{m=1}^5\langle k_m\rangle^{d_m}$ for some
	$\boldsymbol d\in\mathscr I$, since the pairs share at most one
	index. When the two index pairs coincide, both frequencies have exponent two,
	so the resulting monomial lies outside $\mathscr I$. We treat this case
	separately in the estimates below.
	
	The phases $\beta_p$, $\beta_3$, and $\Gamma_p$ are finite linear
	combinations of products $k_i k_q$ with $i<q$. Bounding
	each phase by a constant times the largest pair product and using
	the monotonicity of $x^2/(1+\tau x)$ gives
	\[
	\frac{|\beta|^2}{1+\tau|\beta|}
	\lesssim \sum_{1\le i<q\le 5}\frac{|k_i k_q|^2}{1+\tau|k_i k_q|},
	\qquad \beta\in\{\beta_p,\beta_3\}.
	\]
	For the product $\mu(\tau|\beta_3|)\mu(\tau|\Gamma_p|)$,
	subadditivity reduces the estimate to terms of the form
	$\mu(\tau|k_i k_q|)\mu(\tau|k_u k_v|)$. For distinct pairs,
	$\mu(x)\le x$ gives the corresponding degree-four monomial
	after multiplication by $\tau^2$. For coincident pairs, writing
	$z=|k_i k_q|$, we use
	$$
	\tau^2\mu(\tau z)^2
	\lesssim \frac{\tau^4z^2}{(1+\tau z)^2}
	\le \frac{\tau^4z^2}{1+\tau z}\lesssim \tau^4\sum_{1\le i<q\le 5}\frac{|k_i k_q|^2}{1+\tau|k_i k_q|}.
	$$
	
	For the separable approximation errors, use $a_1=k_1^2$,
	$a_2\le2(k_3^2+k_4^2)$, and $b=k_5^2$ to reduce the product
	$\mu(\tau a_p)\mu(\tau b)$ to terms of the form
	$\mu(\tau k_i^2)\mu(\tau k_5^2)$ with $i\ne5$. For the $F$ error, using $\mu(x)\lesssim x/(1+x)$ gives
	\[
	\tau^2\mu(\tau k_i^2)\mu(\tau k_5^2)
	\lesssim
	\frac{\tau^4|k_i k_5|^2}{(1+\tau k_i^2)(1+\tau k_5^2)}
	\le
	\frac{\tau^4|k_i k_5|^2}{1+\tau|k_i k_5|},
	\]
	where the last inequality follows from
	$(1+x^2)(1+y^2)\ge (1+xy)^2$ for $x,y\ge0$.
	For the $G$ and $Q_p$ errors, expand the prefactor
	$|\beta_p|$ or $|\beta_3|$ into frequency pairs. Using
	$\mu(x)\lesssim\sqrt{x/(1+x)}$ and the same inequality gives
	\[
	\tau^3|k_u k_v|\mu(\tau k_i^2)\mu(\tau k_5^2)
	\lesssim
	\frac{\tau^4|k_u k_v|\,|k_i k_5|}
	{1+\tau|k_i k_5|}.
	\]
	If $\{u,v\}=\{i,5\}$, this gives the required pair bound.
	Otherwise, dropping the denominator gives a monomial in
	$\mathscr I$. The filter errors also give monomials in
	$\mathscr I$, since they are bounded by sums of terms of the form
	$\tau^4|k_i k_q k_u|$ with $i\ne q$.
	
	Combining these bounds yields
	\begin{equation}\label{eq:quintic-defect-kernel-bound}
		|K_{5,p,R}^{j,\tau}|
		\lesssim \tau^4\left(
		\sum_{1\le i<q\le 5}\frac{|k_i k_q|^2}{1+\tau|k_i k_q|}
		+\sum_{\boldsymbol d\in\mathscr I}
		\prod_{\ell=1}^5\langle k_\ell\rangle^{d_\ell}
		\right).
	\end{equation}
	
	Dropping the denominators leaves at most two derivatives on each
	input. For $r>\tfrac52$, the algebra property of $H^{r-2}$ gives
	\[
	\|\mathcal R_{(5)}^\tau(t_n)\|_{H^{r-2}}
	\lesssim
	\tau^4\|v\|_{C([0,T];H^r)}^5.
	\]
	
	For the cumulative estimate, define $\mathcal D_m(f_1,\ldots,f_5)$
	by replacing the five input factors in $V_1$ and $V_2$ of
	$\mathcal R_{(5)}^\tau(t_m)$ with $f_1,\ldots,f_5$, respectively,
	with conjugation retained in the corresponding slots.
	Thus $\mathcal R_{(5)}^\tau(t_m)
	=\mathcal D_m(v(t_m),\ldots,v(t_m))$.
	Set $\mathcal A_n=\sum_{m=0}^{n-1}\mathcal D_m$.
	The contribution of the monomials in $\mathscr I$ to
	$\mathcal A_n$ is bounded by $O_T(\tau^3)$
	from $(H^2)^5$ to $L^2$. Indeed, by Young's
	convolution inequality, we place the most heavily weighted input
	in Fourier $\ell^2$ and the other four in $\ell^1$, using
	$\|\langle k\rangle\widehat f_k\|_{\ell^1}\lesssim\|f\|_{H^2}$.
	
	For each quintic term, let $\Phi=\phi_p$ denote its phase.
	By Lemma~\ref{lem:discrete-time-summation}(i) and
	\eqref{eq:quintic-defect-kernel-bound}, the pair terms in the
	bound for the $k$th Fourier coefficient of
	$\mathcal A_n(f_1,\ldots,f_5)$ are bounded by
	$$
	C_T\tau^3
	\sum_{1\le i<q\le5}
	\sum_{k_1+\cdots+k_5=k}
	\frac{g_\tau(\Phi)|k_i k_q|^2}{1+\tau|k_i k_q|}
	\prod_{a=1}^5 c_a(k_a),
	$$
	where $c_a$ denotes the absolute Fourier coefficients of the
	input in slot $a$, with conjugation as prescribed by $V_p$.
	Fix a pair $(i,q)$ and the other three frequencies, whose
	sum is denoted by $\kappa$. Write $x=k_i$, $y=k_q$, and
	$S=k-\kappa$, so that $x+y=S$. Set $h_a(z)=|z|^2c_a(z)$ for
	$a=i,q$. After taking out the three fixed Fourier factors,
	the remaining sum is
	\[
	H(k)=\sum_{x+y=S}
	\frac{g_\tau(\Phi)}{1+\tau|xy|}\,h_i(x)h_q(y).
	\]
	
	At fixed $k$, the equal-sign case gives $\Phi=C\pm2xy$, which
	is quadratic in $x$ with leading coefficient of modulus $2$.
	The opposite-sign case gives $\Phi=C\pm S(2x-S)$, which is
	affine with slope of modulus $2|S|$. In the equal-sign case,
	$|\Phi-C|=2|xy|$. In the opposite-sign case, since $x+y=S$,
	\[
	\min_{\nu\in\{-1,1\}}|\Phi-(C+\nu S^2)|
	=2|S|\min\{|x|,|y|\}
	\le 4|xy|.
	\]
	Let $\mathcal C=\{C\}$ in the equal-sign case and
	$\mathcal C=\{C-S^2,C+S^2\}$ in the opposite-sign case.
	On each resonance interval $I_\ell$, we have
	$\min_{c\in\mathcal C}|2\pi\ell-\tau c|
	\le\pi+4\tau|xy|$, which implies
	\[
	\frac{1}{(1+\tau|xy|)^2}
	\lesssim
	\sum_{c\in\mathcal C}
	\frac{1}{(1+|2\pi\ell-\tau c|)^2}.
	\]
	
	Except in the opposite-sign case with $S=0$, applying
	Lemma~\ref{lem:discrete-time-summation}(ii) with $\rho=2$
	on each resonance interval, we obtain
	\[
	\sum_{x+y=S}\frac{g_\tau(\Phi)^2}{(1+\tau|xy|)^2}
	\lesssim
	\sum_{c\in\mathcal C}\sum_{\ell\in\mathbb Z}
	\frac{1}{(1+|2\pi\ell-\tau c|)^2}
	\sum_{\substack{x+y=S\\\Phi\in I_\ell}}g_\tau(\Phi)^2
	\lesssim
	\sum_{c\in\mathcal C}\sum_{\ell\in\mathbb Z}
	\frac{1}{(1+|2\pi\ell-\tau c|)^2}
	\lesssim1.
	\]
	
	Since $S=k-\kappa$, the condition $k\ne\kappa$ is equivalent to
	$S\ne0$. For these outputs, Cauchy--Schwarz at fixed $k$ gives
	$
	|H(k)|^2\lesssim
	\sum_{x+y=k-\kappa}h_i(x)^2h_q(y)^2.
	$
	Summing over $k\ne\kappa$ yields
	$
	\sum_{k\ne\kappa}|H(k)|^2
	\lesssim \|h_i\|_{\ell^2}^2\|h_q\|_{\ell^2}^2,
	$
	since each pair $(x,y)$ determines a unique output
	$k=x+y+\kappa$. At $k=\kappa$, we have $S=0$ and hence $y=-x$.
	Using $g_\tau\le1$ and Cauchy--Schwarz, we obtain
	$
	|H(\kappa)|\le\sum_{x\in\mathbb Z}h_i(x)h_q(-x)
	\le\|h_i\|_{\ell^2}\|h_q\|_{\ell^2}.
	$
	Combining the two bounds gives
	$\|H\|_{\ell^2_k}\lesssim
	\|h_i\|_{\ell^2}\|h_q\|_{\ell^2}$ uniformly in the other
	three frequencies.
	
	Let $J=\{a_1,a_2,a_3\}=\{1,\ldots,5\}\setminus\{i,q\}$,
	with $a_1<a_2<a_3$, and set
	$\boldsymbol{k}_J=(k_{a_1},k_{a_2},k_{a_3})\in\mathbb Z^3$.
	Write $H_{\boldsymbol{k}_J}$ for the sequence $H$ with these
	three frequencies fixed.
	By the triangle inequality in $\ell^2_k$ and the uniform bound
	for $H_{\boldsymbol k_J}$, we have
	\[
	\begin{aligned}
		\left\|
		\sum_{\boldsymbol k_J\in\mathbb Z^3}
		\left(\prod_{a\in J}c_a(k_a)\right)H_{\boldsymbol k_J}
		\right\|_{\ell^2_k}
		\le
		\sum_{\boldsymbol k_J\in\mathbb Z^3}
		\left(\prod_{a\in J}c_a(k_a)\right)
		\|H_{\boldsymbol k_J}\|_{\ell^2_k}\lesssim
		\|h_i\|_{\ell^2}\|h_q\|_{\ell^2}
		\prod_{a\in J}\|c_a\|_{\ell^1}.
	\end{aligned}
	\]
	Since $\|h_b\|_{\ell^2}\lesssim\|f_b\|_{H^2}$ for
	$b\in\{i,q\}$ and $\|c_a\|_{\ell^1}\lesssim\|f_a\|_{H^2}$
	for $a\in J$, restoring the factor $C_T\tau^3$ bounds the
	contribution of the fixed pair by
	$$
	C_T\tau^3\|h_i\|_{\ell^2}\|h_q\|_{\ell^2}
	\prod_{a\in J}\|c_a\|_{\ell^1}
	\lesssim_T \tau^3\prod_{a=1}^5\|f_a\|_{H^2}.
	$$
	Summing over $1\le i<q\le5$ and adding the contribution of
	the monomials indexed by $\mathscr I$, we obtain
	\[
	\max_{0\le n\le L}\|\mathcal A_n(f_1,\ldots,f_5)\|_{L^2}
	\lesssim_T \tau^3\prod_{a=1}^5\|f_a\|_{H^2}.
	\]
	
	Duhamel's formula for the twisted variable directly yields
	$$
	\sum_{m=0}^{L-2}\|v(t_{m+1})-v(t_m)\|_{H^2}
	\lesssim T\|v\|_{C([0,T];H^2)}^3.
	$$
	The cumulative estimate follows from
	Lemma~\ref{lem:discrete-time-summation}(iii) with
	$d=5$, $X=H^2$, $Y=L^2$, and $w_m=v(t_m)$.
\end{proof}
\subsubsection{Proof of Lemma~\ref{lemma3.4}}\label{pooflemma3.4}
\begin{proof}
	Let $\theta_j$ denote the retained phase in the $j$th row of
the preceding table. Since the remaining oscillatory factors
$\mathfrak a(\sigma)$ are evaluated at $\sigma=0$ in the
approximation, each remainder kernel has the form
	\[
	K_{7,R}^{j,\tau}
	=\int_0^\tau \sigma^2 \mathrm e^{\ii\sigma\theta_j}
	\bigl(\mathfrak a(\sigma)-\mathfrak a(0)\bigr)\,\mathrm d\sigma.
	\]
	For the first two summands, $\mathfrak a(\sigma)=\mathrm  e^{2\ii\sigma\vartheta}$,
	where $\vartheta$ is the corresponding quadratic residual phase
	and $\vartheta=\beta_4$ for the first summand.
For the last two, expanding the factors involving $\Omega_1^n$
gives amplitudes of the form
\[
\mathfrak a(\sigma)=\mathrm e^{2\ii\sigma\vartheta}
\chi_\mu(2\ii\sigma\eta k_3^2)\chi_\nu(2\ii\sigma k_6^2),
\qquad
\chi_0=\varphi_1,\quad \chi_1=\psi_1,
\quad \mu,\nu\in\{0,1\},
\]
where $\eta=-1$ for the third summand and $\eta=1$ for the fourth.
	Using $|\mathrm e^{\ii y}-1|\lesssim |y|$, the boundedness of $\chi_b$
	on the imaginary axis, and
	$|\chi_b(\ii y)-\chi_b(0)|\lesssim |y|$ for $b=0,1$,
	we obtain in all four cases
	\[
	|\mathfrak a(\sigma)-\mathfrak a(0)|
	\lesssim \sigma\bigl(|\vartheta|+k_3^2+k_6^2\bigr)
	\lesssim \sigma\sum_{a=1}^7\langle k_a\rangle^2.
	\]
	Consequently, summing over the individual remainder kernels gives
	$|K_{7,R}^\tau|\lesssim
	\tau^4\sum_{a=1}^7\langle k_a\rangle^2$.
	By Lemma~\ref{kp}~\textup{(ii)} with $\delta=r-2$ and
	$\gamma=2$, followed by \textup{(i)} with $\gamma=r$, we obtain
	\[
	\begin{aligned}
		\|\mathcal R_{(7)}^\tau(t_n)\|_{H^{r-2}}\lesssim \tau^4\|v\|_{C([0,T];H^r)}^7.
	\end{aligned}
	\]
	The $L^2$ estimate follows from $\delta=r-2=0$.
\end{proof}
\subsubsection{Proof of Lemma~\ref{lemma3.1}}\label{pooflemma3.1}
\begin{proof}
	A comparison of \eqref{integral} and \eqref{expressions}
	shows that the remainder $\mathcal{R}_{(0)}^\tau(t_n)$ can be
	written as a sum of terms of the following form with indices
	$j,k,l\in\{0,1,2,3\}$ satisfying $j+k+l\ge3$:
	\begin{equation}\label{eq:4.-1}
		-\ii\int_0^\tau
		S_{-(t_n+\sigma)}
		\Big[
		S_{-(t_n+\sigma)}
		\overline{\mathcal W_j}(t_n)\cdot
		S_{t_n+\sigma}
		\mathcal W_k(t_n)\cdot
		S_{t_n+\sigma}
		\mathcal W_l(t_n)
		\Big]\mathrm d\sigma,
	\end{equation}
	where
	$\mathcal W_0(t_n)=v(t_n)$,
	$\mathcal W_1(t_n)=\Omega_1^n(v(t_n),\sigma)$,
	$\mathcal W_2(t_n)=\Omega_2^n(v(t_n),\sigma)$, and
	$\mathcal W_3(t_n)=\mathcal R_{0,0}^\sigma(t_n)$.
	
	Since $\varphi_1$ and $\psi_1$ are bounded on the imaginary
	axis, applying Lemma~\ref{kp}~\textup{(i)} for $r\ge2$ yields
	\begin{equation}\label{eq:4.0}
		\begin{aligned}
			\|\mathcal W_1(t_n)\|_{H^r}
			&=\|\Omega_1^n(v(t_n),\sigma)\|_{H^r}
			\lesssim\sigma\|v\|_{C([0,T];H^r)}^3,\\
			\|\mathcal W_2(t_n)\|_{H^r}
			&=\|\Omega_2^n(v(t_n),\sigma)\|_{H^r}
			\lesssim\sigma^2\|v\|_{C([0,T];H^r)}^5.
		\end{aligned}
	\end{equation}
	
	We first establish the one-step error bound in $H^{r-2}$
	for $r>\tfrac52$. Combining \eqref{eq:R00}, \eqref{eq:4.-1}, 
	and \eqref{eq:4.0} with Lemma~\ref{kp}~\textup{(i)} and
	using $j+k+l\ge3$, we obtain
	\[
	\|\mathcal{R}_{(0)}^\tau(t_n)\|_{H^{r-2}}
	\le C_r\tau^4.
	\]
	
	We now establish the cumulative bound for
	$v\in C([0,T];H^2)$. Combining \eqref{eq:3.1}, \eqref{ex},
	and \eqref{eq:4.0} with $r=2$, we have
	$
	\|\mathcal W_3(t_n)\|_{H^2}
	\le C\bigl(\|v\|_{C([0,T];H^2)}\bigr)\sigma.
	$
	In what follows, $\mathcal O_{L^2}(s^p)$ denotes a remainder
	whose $L^2$ norm is bounded by
	$C(\|v\|_{C([0,T];H^2)})s^p$, uniformly in $n$ and $s$.
	
	For brevity, set $U_s=S_{t_n+s}v(t_n)$.
	Freezing the nonlinearity at $t_n$ and expanding the cubic
	product with the aid of Lemma~\ref{kp}~\textup{(i)} and
	$\|(\mathrm e^{\ii s\partial_x^2}-1)f\|_{L^2}
	\le s\|f\|_{H^2}$ yields
	\[
	|u(t_n+s)|^2u(t_n+s)-|U_s|^2U_s
	=-\ii s|U_0|^4U_0+\mathcal O_{L^2}(s^2).
	\]
	Substituting this expansion into \eqref{eq:3.1} with $\tau$
	replaced by $\sigma$ and using \eqref{Omega2} together
	with the same group increment estimate, we obtain
	\[
	\mathcal W_3(t_n)
	=\Big(-\ii\int_0^\sigma
		S_{-(t_n+s)}
	\bigl(|U_s|^2U_s\bigr)\,\mathrm ds
	-\Omega_1^n(v(t_n),\sigma)\Big)+\mathcal O_{L^2}(\sigma^3)=:\mathcal C^\sigma(t_n)+\mathcal O_{L^2}(\sigma^3).
	\]
	Following the argument in \cite{Ostermann2022}, we write
	\[
	\begin{aligned}
		\mathcal C^\sigma(t_n)
		=-\ii\sum_{k\in\mathbb Z}
		\sum_{k_1+k_2+k_3=k}
		\mathrm e^{\ii t_n\phi}
		\{\hat{\bar v}_{k_1}\hat v_{k_2}\hat v_{k_3}\}(t_n)
		\mathrm e^{\ii kx}\cdot\int_0^\sigma
		\mathrm e^{2\ii s k_1^2}
		\left[
		\mathrm e^{\ii s\beta}-1
		-\frac{s}{\sigma}
		\bigl(\mathrm e^{\ii\sigma\beta}-1\bigr)
		\right]\mathrm ds,
	\end{aligned}
	\]
	where $\beta=2(k_1k_2+k_1k_3+k_2k_3)$.
	By Lemma~\ref{jb} and the same frequency-product estimate
	as in the quintic case, we readily obtain
	\[
	\left|\int_0^\sigma
	\mathrm e^{2\ii s k_1^2}
	\left[
	\mathrm e^{\ii s\beta}-1
	-\frac{s}{\sigma}
	\bigl(\mathrm e^{\ii\sigma\beta}-1\bigr)
	\right]\mathrm ds\right|
	\lesssim\frac{\sigma^3|\beta|^2}{1+\sigma|\beta|}
	\lesssim \sigma^2\sum_{1\le a<b\le3}|k_a k_b|.
	\]
	Applying Lemma~\ref{kp}~\textup{(i)} then yields
	$\mathcal W_3(t_n)=\mathcal O_{L^2}(\sigma^2)$.
	
	Returning to \eqref{eq:4.-1}, if $j,k,l\le2$, the algebra
	property of $H^2$ bounds the integrand in $L^2$ by
	$C\sigma^{j+k+l}\le C\sigma^3$. If one index equals $3$
	and another is nonzero, placing one factor $\mathcal W_3$
	in $L^2$ and the remaining factors in $L^\infty$ also
	gives an $\mathcal O_{L^2}(\sigma^3)$ bound.
	
	It remains to consider the triples $(3,0,0)$, $(0,3,0)$,
	and $(0,0,3)$. Replacing $\mathcal W_3(t_n)$ by
	$\mathcal C^\sigma(t_n)$ in these three terms, we write
	$
	\mathcal{R}_{(0)}^\tau(t_n)
	=\mathcal R_{(0)}^{1,\tau}(t_n)+\mathcal R_{(0)}^{2,\tau}(t_n),
	$
	where
	\[
	\begin{aligned}
		\mathcal R_{(0)}^{1,\tau}(t_n)
		=-\ii\int_0^\tau
		S_{-(t_n+\sigma)}
		\Big[
		U_\sigma^2
		S_{-(t_n+\sigma)}
		\overline{\mathcal C^\sigma(t_n)}+2|U_\sigma|^2
		S_{t_n+\sigma}
		\mathcal C^\sigma(t_n)
		\Big]\mathrm d\sigma,
	\end{aligned}
	\]
	and the preceding bounds give
	$
	\mathcal R_{(0)}^{2,\tau}(t_n)=\mathcal{O}_{L^2}(\tau^4).
	$
	
	The intermediate frequency square cancels between the inner
	and outer phases at time $t_n$, so the resulting phase has
	the same quadratic structure as in the quintic case.
	Moreover, the corresponding Fourier multipliers are bounded by
	\[
	\int_0^\tau
	\frac{\sigma^3|\beta|^2}{1+\sigma|\beta|}
	\,\mathrm d\sigma
	\lesssim
	\frac{\tau^4|\beta|^2}{1+\tau|\beta|}
	\lesssim
	\tau^4\sum_{1\le i<\ell\le5}
	\frac{|k_i k_\ell|^2}{1+\tau|k_i k_\ell|}.
	\]
	The fixed-input partial-sum estimate established in the
	quintic case, together with
	Lemma~\ref{lem:discrete-time-summation}~\textup{(iii)}
	and the $H^2$ increment bound following from \eqref{eq:3.1},
	yields
	$
	\max_{1\le n\le L}
	\left\|\sum_{m=0}^{n-1}\mathcal R_{(0)}^{1,\tau}(t_m)\right\|_{L^2}
	\le C_T\tau^3.
	$
	Since $\mathcal R_{(0)}^{2,\tau}(t_m)=\mathcal O_{L^2}(\tau^4)$,
	we conclude that
	\[
	\max_{1\le n\le L}
	\left\|\sum_{m=0}^{n-1}\mathcal R_{(0)}^\tau(t_m)\right\|_{L^2}
	\le C_T\tau^3,
	\]
	where $C_T$ depends only on $T$ and $\|v\|_{C([0,T];H^2)}$.
\end{proof}
\subsection{Proof of the error estimate~\eqref{T1.1} in $H^{r}$ with $r>\tfrac{1}{2}$}\label{r12}
By taking the difference between the numerical scheme~\eqref{vscheme} and the exact solution, we arrive at
\begin{equation} \label{error expression}
v^{n+1} - v(t_{n+1}) 
= \mathcal{L}^n + \Psi^n(v^n,\tau) - \Psi^n \big(v(t_n),\tau\big),
\end{equation}
where $\mathcal{L}^n = \Psi^n\big(v(t_n),\tau \big) - v(t_{n+1})$ is the local error. For convenience, throughout the remainder of Section~\ref{proof} we set $g_n := v^n - v(t_n)$.
\subsubsection{Local error estimate in $H^r$ with $r>\frac12$}\label{4.2.1}
We derive the local error bounds for \eqref{vscheme}, as stated in the lemma below.
\begin{lemma}[Local error]\label{lem:local-error}
	Let $0<\tau\le1$, $r>\tfrac12$, and
	$v\in C([0,T];H^{r+2})$. Then
	\[
	\|\mathcal L^n\|_{H^r}\le C\tau^4,
	\qquad 0\le n<L,
	\]
	where the constant $C$ depends only on $r$ and
	$\|v\|_{C([0,T];H^{r+2})}$.
\end{lemma}
\begin{proof}
It follows easily from \eqref{expression} that
$
\mathcal{L}^n=-\sum_{j\in\{0,3,5,7\}}\mathcal R_{(j)}^\tau(t_n).$
Hence, the desired estimates are direct consequences of Lemmas~\ref{lemma3.1}--\ref{lemma3.4}.
\end{proof}
\subsubsection{Stability}
Next, we establish the stability of the numerical scheme \eqref{vscheme}.
It follows directly from \eqref{expression} that
\begin{equation}\label{psiexpression}
\|\Psi^n(v(t_n),\tau)-\Psi^n(v^n,\tau)\|_{H^r}\le \| g_n \|_{H^r}+\sum\limits_{j\in\{3,5,7\}}	\|\Psi_{(j)}^n(v^n,\tau)-\Psi_{(j)}^n(v(t_n),\tau)\|_{H^r}.
\end{equation}
 Due to the distinct structure of $\Psi_{(3)}^{n}$, we estimate it separately from the remaining terms
 $\sum_{j\in\{5,7\}}\Psi_{(j)}^{n}$. Combining these two estimates yields the desired stability bound in
 $H^r$ ($r > \frac{1}{2}$) for the scheme \eqref{vscheme}.
 In what follows, $C>0$ denotes a generic constant whose value may change from line to line,
 but depends only on $\|v\|_{C([0,T];H^r)}$.
 We first establish the following stability estimate for $\Psi_{(3)}^{n}$ in $H^r$.
\begin{lemma}\label{r:psi1ns}
Let $r>\tfrac{1}{2}$, $0<\tau\le1$ and $v\in C([0,T];H^r)$. We have  
\begin{align*}
\|\Psi_{(3)}^n(v^n,\tau)-\Psi_{(3)}^n(v(t_n),\tau)\|_{H^r}\le C\tau \Big(\| g_n\|_{H^r}+\|  g_n\|_{H^r}^3\Big),
\end{align*}
where the constant $C$ depends only on $\|v\|_{C([0,T];H^r)}$.
\end{lemma}
\begin{proof}
	Using the Fourier representation in \cite{cao2024new}, we write
	\begin{equation}\label{psiall}
		\Psi_{(3)}^n(v^n,\tau)-\Psi_{(3)}^n(v(t_n),\tau)
		=\sum_{j=1}^5\Psi_{(3)}^{j,n}.
	\end{equation}
	The first term corresponds to the zero Fourier mode,
	whose multiplier is bounded by $\tau$. For $k\ne0$,
	after factoring out $\mathrm e^{\ii t_n\phi}$,
	the remaining multipliers are given by
	\begin{align*}
		G_2&=-2\int_0^\tau\frac{k_1+k_2}{k}
		\mathrm e^{\ii s\phi}\,\mathrm ds,
		&
		G_3&=\int_0^\tau\frac{k_1}{k}
		\mathrm e^{2\ii skk_1}\,\mathrm ds,\\
		G_4&=\int_0^\tau\frac{k_1}{k}
		\bigl(\mathrm e^{2\ii sk_2k_3}-1\bigr)\,\mathrm ds,
		&
		G_5&=\frac{k_1}{k}\mathcal M_\tau(kk_1)
		\bigl(\mathrm e^{2\ii\tau k_2k_3}-1\bigr).
	\end{align*}
	Since $|\mathcal M_\tau(kk_1)|\lesssim\tau$,
	all five Fourier multipliers are bounded by
	\begin{equation}\label{5}
		C\tau\frac{\sum_{a=1}^3\langle k_a\rangle}{\langle k\rangle}.
	\end{equation}
	
	Following \cite{li2021fully}, we expand the difference of
	the cubic products into terms of the form $f_1f_2f_3$,
	where one factor is $g_n$ or its conjugate and the remaining
	two are $v^n$, $v(t_n)$, or their conjugates.
	Replacing their Fourier coefficients by absolute values,
	and using \eqref{5} together with
	Lemma~\ref{kp}~\textup{(i)} and Lemma~\ref{li}~\textup{(i)},
	we bound each such contribution in $H^r$ by
	\[
	\tau\sum_{a=1}^3
	\left\|J^{-1}\left(Jf_a\prod_{b\ne a}f_b\right)\right\|_{H^r}
	\lesssim\tau\prod_{a=1}^3\|f_a\|_{H^r}.
	\]
	Consequently, by \eqref{psiall}, we obtain
	\[
	\begin{aligned}
		\|\Psi_{(3)}^n(v^n,\tau)-\Psi_{(3)}^n(v(t_n),\tau)\|_{H^r}\le C\tau\bigl(\|g_n\|_{H^r}+\|g_n\|_{H^r}^3\bigr).
	\end{aligned}
	\]
\end{proof}
Having established the stability estimate for $\Psi_{(3)}^{n}$, we next derive an $H^r$ stability estimate for the remaining terms.
\begin{lemma}\label{r:rs}
Suppose that $r>\tfrac{1}{2}$, $0<\tau\le1$ and $v\in C([0,T];H^r)$. It follows that
$$\sum\limits_{j\in\{5,7\}}\|\Psi_{(j)}^n(v^n,\tau)-\Psi_{(j)}^n(v(t_n),\tau)\|_{H^r}\le C\tau \Big(\| g_n\|_{H^r}+\|  g_n\|_{H^r}^7\Big),$$
where the constant $C$ depends only on $\|v\|_{C([0,T];H^r)}$.
\end{lemma}
\begin{proof}
	Utilizing the uniform boundedness of the Fourier multipliers appearing in $\Psi_{(5)}^n$ and $\Psi_{(7)}^n$ on $H^r$, the fact that $S_t$ is an isometry on $H^r$, and Lemma~\ref{kp}\textup{(i)}, we infer from the definition of $\Psi_{(7)}^n$ that
	$
	\|\Psi_{(7)}^n(v^n,\tau)-\Psi_{(7)}^n(v(t_n),\tau)
	\|_{H^r}
	\le C\tau^3\Big(\|g_n\|_{H^r}+\|g_n\|_{H^r}^7\Big).
	$
	For $\Psi_{(5)}^n$, we additionally use the Leibniz rule and the filtered estimates
	$
	\|P_\tau f\|_{H^r}\le\|f\|_{H^r},\,
	\|DP_\tau f\|_{H^r}\lesssim\tau^{-1}\|f\|_{H^r}.
	$
	Since each term with the prefactor $\tau^3$ contains two spatial derivatives acting on the filtered factors after applying the Leibniz rule, the preceding arguments yield
	$
	\|\Psi_{(5)}^n(v^n,\tau)-\Psi_{(5)}^n(v(t_n),\tau)\|_{H^r}
	\le C\tau\Big(\|g_n\|_{H^r}+\|g_n\|_{H^r}^5\Big).
	$
	Combining these estimates and using
	$\|g_n\|_{H^r}^5\le \|g_n\|_{H^r}+\|g_n\|_{H^r}^7$
	gives the desired result.
\end{proof}
\par Combining the two estimates derived above, we obtain the following stability estimate for \eqref{vscheme}.
\begin{lemma}[Stability]\label{stability}
Letting $r>\tfrac{1}{2}$ and $0<\tau\le1$, one has
$$\|\Psi^n(v^n,\tau)-\Psi^n(v(t_n),\tau)\|_{H^r}\le(1+C\tau)\| g_n\|_{H^r}+C\tau \| g_n\|_{H^r}^7,$$
where the constant $C$ depends only on $\| v \|_{C([0,T];H^r)}$.
\end{lemma}
\begin{proof}
The desired result follows directly from \eqref{psiexpression} and Lemmas~\ref{r:psi1ns}-\ref{r:rs}.
\end{proof}
\subsubsection{Proof of the statement \eqref{T1.1}}\label{proofof1.4}
\begin{proof}For $0<\tau\le 1$, combining Lemma~\ref{lem:local-error} and Lemma~\ref{stability}, \eqref{error expression} in fact shows that there exists a constant $C$ depending only on $\|v\|_{C([0,T];H^{r+2})}$ such that $$
\|v(t_{n+1})-v^{n+1}  \|_{H^r}\le C\tau^4+(1+C\tau)\|g_n\|_{H^r}+C\tau\|g_n\|_{H^r}^7,\  0\le n< L.$$
Iterating the recursion yields the bound $
\| v(t_{n+1})-v^{n+1}\|_{H^r}\le C\tau \sum_{j=0}^n (1+C\tau)^j \Big[\| g_{n-j} \|_{H^r}^7 +\tau^3 \Big].$
From the above estimate, we deduce, by a bootstrap argument, that there exist $\tau_0>0$ and $C_0>0$ such that for any $0<\tau\le \tau_0$:
$$
\|v(t_{n+1})-v^{n+1}  \|_{H^r}\le C\tau^4 \sum_{j=0}^n (1+C\tau )^j\le C_0\tau^3,\ n\ge0.$$
Observe that the constants $\tau_0$ and $C_0$ depend only on $T$ and $\|v\|_{C([0,T];H^{r+2})}$. Moreover, since $S_t$ is a linear isometry on $H^s$ for $s\ge 0$ and $t\in\mathbb{R}$, the statement~\eqref{T1.1} follows immediately.\end{proof}
\subsection{Proof of the error estimate~\eqref{T1.2} in $L^2$}\label{r13}
This subsection proves the global convergence result~\eqref{T1.2} in $L^2$.
We first establish a cumulative local error bound in $L^2$.
Building on Lemma~\ref{stability}, we then derive an a priori $H^1$ bound
for the numerical solution, which is used to obtain an $L^2$ stability estimate.
The convergence result~\eqref{T1.2} follows by combining the cumulative
local error bound with the $L^2$ stability estimate.
\subsubsection{Cumulative local error estimate in $L^2$}
We derive the cumulative local error bound for \eqref{vscheme} in $L^2$, as stated in the lemma below.
\begin{lemma}[Cumulative local error]\label{lem:cumulative-local-error}
	Let $0<\tau\le1$ and
	$v\in C([0,T];H^{2})$. Then
	\[
	\max_{1\le n\le L}
	\left\|\sum_{m=0}^{n-1}\mathcal L^m\right\|_{L^2}
	\le C\tau^3,
	\]
	where the constant $C$ depends only on $T$ and
	$\|v\|_{C([0,T];H^{2})}$.
\end{lemma}
\begin{proof}
	It follows from \eqref{expression} that
	$
	\sum_{m=0}^{n-1}\mathcal L^m
	=-\sum_{j\in\{0,3,5,7\}}\sum_{m=0}^{n-1}
	\mathcal R_{(j)}^\tau(t_m).
	$
	By Lemmas~\ref{lemma3.1}--\ref{lem:quintic-remainder}, we have
	\[
	\max_{1\le n\le L}
	\left\|\sum_{m=0}^{n-1}\mathcal R_{(j)}^\tau(t_m)\right\|_{L^2}
	\le C\tau^3,\qquad j\in\{0,3,5\}.
	\]
	Moreover, applying Lemma~\ref{lemma3.4} with $r=2$ yields
	\[
	\max_{1\le n\le L}
	\left\|\sum_{m=0}^{n-1}\mathcal R_{(7)}^\tau(t_m)\right\|_{L^2}
	\le
	\sum_{m=0}^{L-1}\|\mathcal R_{(7)}^\tau(t_m)\|_{L^2}
	\lesssim T\tau^3\|v\|_{C([0,T];H^2)}^7.
	\]
	Combining these estimates gives the desired result.
\end{proof}

\subsubsection{A priori $H^1$ bound for the numerical solution}
To derive an a priori $H^1$ bound for the numerical solution, we first establish a local error estimate in $H^1$ under a weaker regularity assumption $u^0\in H^{2}$. This is stated in the following lemma, with $\mathcal{L}^n$ defined in Section~\ref{r12}.
\begin{lemma}\label{lem:local-H1}
Assuming $0<\tau\le 1$ and $v\in C([0,T];H^{2})$, then it is obtained that
$
\|\mathcal{L}^n\|_{H^1}\le C\tau^2,$
where the constant $C$ depends only on 
$\|v\|_{C([0,T];H^{2})}$.
\end{lemma}
\begin{proof}
	By \eqref{expression}, it suffices to estimate
	$\mathcal R_{(j)}^\tau(t_n)$ in $H^1$ for
	$j\in\{0,3,5,7\}$.
	We first consider the cubic remainder. Using
	\eqref{eq:cubic-remainder-kernel-bound}, we obtain for $k\ne0$
	\[
	\begin{aligned}
		\langle k\rangle |K_{3,R}^\tau|
		\lesssim
		\frac{\tau^4|kk_1||k_2k_3||k_1k_2k_3|}
		{(1+\tau|kk_1|)(1+\tau^2|k_2k_3|^2)}\lesssim\tau^2|k_1k_2k_3|.
	\end{aligned}
	\]
Since the zero mode vanishes,
	Lemma~\ref{kp}\textup{(i)} yields
	$
	\|\mathcal R_{(3)}^\tau(t_n)\|_{H^1}
	\lesssim\tau^2\|v(t_n)\|_{H^2}^3.
	$
	
	For the quintic and septic remainders, combining
	\eqref{eq:4.0} with Lemma~\ref{kp}\textup{(i)} and using
	the arguments in the proof of Lemma~\ref{r:rs} with
	$\|DP_\tau f\|_{H^1}\le\|f\|_{H^2}$, we obtain
	\[
	\begin{aligned}
		\|\mathcal R_{(5)}^\tau(t_n)\|_{H^1}
		&\le\|I_{(5)}^\tau(t_n)\|_{H^1}
		+\|\Psi_{(5)}^n(v(t_n),\tau)\|_{H^1}
		\lesssim\tau^2\|v\|_{C([0,T];H^2)}^5,\\
		\|\mathcal R_{(7)}^\tau(t_n)\|_{H^1}
		&\le\|I_{(7)}^\tau(t_n)\|_{H^1}
		+\|\Psi_{(7)}^n(v(t_n),\tau)\|_{H^1}
		\lesssim\tau^3\|v\|_{C([0,T];H^2)}^7.
	\end{aligned}
	\]
	
	Finally, by \eqref{eq:3.1}, \eqref{ex} and \eqref{eq:4.0},
	we have
	$
	\|\mathcal W_3(t_n)\|_{H^1}
	\le C\big(\|v\|_{C([0,T];H^2)}\big)\sigma.
	$
	Combining this estimate with \eqref{eq:4.-1},
	\eqref{eq:4.0} and Lemma~\ref{kp}\textup{(i)},
	and noting that $j+k+l\ge3$, we obtain
	\[
	\|\mathcal R_{(0)}^\tau(t_n)\|_{H^1}
	\le C\big(\|v\|_{C([0,T];H^2)}\big)\tau^2.
	\]
	Combining the above estimates gives the desired result.
\end{proof}

We proceed by deriving an error bound for the numerical solution in $H^1$, and subsequently use this bound to show that the numerical solution is uniformly bounded in $H^1$ with respect to $\tau$.
\begin{lemma}\label{h1bound}
Let $u^0 \in H^{2}$ and $u^n, n=0,1,\ldots,L$ be the numerical solution of the NLS equation given by \eqref{littlescheme}--\eqref{numericalscheme}. Then there exist positive constants $\tau_1$ and $C_1$ such that for any 
$\tau \in (0,\tau_1]:$ $$
\max_{1 \leq n \leq L} \| u(t_n,\cdot) - u^n \|_{H^1} \leq C_1 \tau,$$
where $C_1$ and $\tau_1$ depend only on the final time $T$ 
and $\| u\|_{C([0,T];H^{2})}$.
\end{lemma}
\begin{proof}
Combining Lemmas~\ref{stability} and \ref{lem:local-H1}, the desired result  follows by the same discrete Gronwall (Lady Windermere's fan) argument as in Section~\ref{proofof1.4}.
\end{proof}

By the triangle inequality and the boundedness of the exact solution in $H^1$, it follows that
$
\|v^n\|_{H^1}
\le \|v(t_n)-v^n\|_{H^1} + \|v(t_n)\|_{H^1}
\lesssim 1$
uniformly for $0\le n\le L$.
This observation can be summarized in the following lemma.
\begin{lemma}\label{1bound}
Let $u^0 \in H^{2}$ and $u^n, n=0,1,\ldots,L$ be the numerical solution of the NLS equation given by \eqref{littlescheme}--\eqref{numericalscheme}. Then there exists a positive constant $\tau_1$ such that for any 
$\tau \in (0,\tau_1]$ the following estimate holds $$
\max_{0\le n\le L}\|u^{n}\|_{H^{1}} \lesssim 1.$$
\end{lemma}
\subsubsection{Stability}
Using the $H^1$ boundedness of the numerical solution established above, 
we derive an $L^2$ stability estimate for the scheme. 
We first treat $\Psi_{(3)}^n$ separately and establish its $L^2$ stability bound, 
then derive corresponding bounds for $\Psi_{(5)}^n$ and $\Psi_{(7)}^n$. 
Combining these yields the $L^2$ stability of the full scheme. 
We begin with the following estimate for $\Psi_{(3)}^{n}$ in $L^2$.
\begin{lemma}\label{0psi1}
Let $v\in C([0,T];H^2)$. Then there exist constants $\tau_1\in(0,1]$ and $C$ such that for any $\tau\in(0,\tau_1]$:
$$
\|\Psi_{(3)}^n(v^n,\tau)-\Psi_{(3)}^n(v(t_n),\tau)\|_{L^2}\le C\tau \| g_n\|_{L^2},$$
where the constants $\tau_1$ and $C$ depend only on $T$ and $\|v\|_{C([0,T];H^2)}$.
\end{lemma}
\begin{proof}
	Using the multiplier bound established in the proof of
	Lemma~\ref{r:psi1ns}, we bound each contribution in $L^2$ by
	\[
	\tau\sum_{a=1}^3
	\left\|J^{-1}\left(Jf_a\prod_{b\ne a}f_b\right)\right\|_{L^2}.
	\]
	By Lemma~\ref{li}~\textup{(ii)}, we have
	\[
	\left\|J^{-1}\left(Jf_a\prod_{b\ne a}f_b\right)\right\|_{L^2}
	\lesssim
	\min\left\{
	\|f_a\|_{L^2}\left\|\prod_{b\ne a}f_b\right\|_{H^1},
	\|f_a\|_{H^1}\left\|\prod_{b\ne a}f_b\right\|_{L^2}
	\right\}.
	\]
	Applying Lemma~\ref{kp}~\textup{(i)} and \textup{(ii)},
	we place the $L^2$ norm on the factor containing $g_n$
	and the $H^1$ norms on the remaining two factors.
	Consequently, by \eqref{psiall} and the uniform $H^1$
	boundedness of the numerical solution, we obtain
	\begin{align*}
		\|\Psi_{(3)}^n(v^n,\tau)-\Psi_{(3)}^n(v(t_n),\tau)\|_{L^2}\le C\tau\|g_n\|_{L^2}.
	\end{align*}
	This is the desired result.
\end{proof}
\par We next establish a stability estimate for the remaining terms in $L^2$.
\begin{lemma}\label{0psi2-7}
Let $v\in C([0,T];H^2)$. Then there exist constants $\tau_1\in(0,1]$ and $C$ such that for any $\tau\in(0,\tau_1]$:
$$
\sum\limits_{j\in\{5,7\}}\|\Psi_{(j)}^n(v^n,\tau)-\Psi_{(j)}^n(v(t_n),\tau)\|_{L^2}\le C\tau \| g_n\|_{L^2},$$
where the constants $\tau_1$ and $C$ depend only on $T$ and $\|v\|_{C([0,T];H^2)}$.
\end{lemma}
\begin{proof}
	By following the proof of Lemma~\ref{r:rs}, with Lemma~\ref{kp}\textup{(i)} replaced by Lemma~\ref{kp}\textup{(ii)}, and always placing the $L^2$ norm on $g_n$ or its complex conjugate, we obtain
	\begin{align*}\label{0psi2n}
		\|\Psi_{(7)}^n(v^n,\tau)-\Psi_{(7)}^n(v(t_n),\tau)\|_{L^2}
		\lesssim\tau^3\|g_n\|_{L^2}
		\big(\|v^n\|_{H^1}^6+\|v(t_n)\|_{H^1}^6\big)
		\lesssim\tau^3\|g_n\|_{L^2}.
	\end{align*}
	For $\Psi_{(5)}^n$, we further employ the filtered derivative estimate
	$
	\|DP_\tau f\|_{L^2}
	\lesssim\tau^{-1}\|f\|_{L^2},
	$
	which, together with the preceding arguments, yields
	$$
	\|\Psi_{(5)}^n(v^n,\tau)-\Psi_{(5)}^n(v(t_n),\tau)\|_{L^2}
	\lesssim\tau\|g_n\|_{L^2}
	\big(\|v^n\|_{H^1}^4+\|v(t_n)\|_{H^1}^4\big)
	\lesssim\tau\|g_n\|_{L^2}.
	$$
	Here, we have used the uniform $H^1$ boundedness of the numerical solution established in Lemma~\ref{1bound}.
	Combining these estimates gives the desired result.
\end{proof}
\begin{lemma}\label{0stability}
For $v\in C([0,T];H^2)$, there exist constants $\tau_1\in(0,1]$ and $C$ such that for any $\tau\in(0,\tau_1]$, the result
$$
\|\Psi^n(v^n,\tau)-\Psi^n(v(t_n),\tau)\|_{L^2}\le(1+C\tau)\| g_n\|_{L^2}$$ is true,
where the constants $\tau_1$ and $C$ only depend on $T$ and $\| v \|_{C([0,T];H^2)}$.
\end{lemma}
\begin{proof}
From \eqref{expression}, it follows that  $$
\|\Psi^n(v(t_n),\tau)-\Psi^n(v^n,\tau)\|_{L^2}\le \| g_n \|_{L^2}+\sum_{j\in\{3,5,7\}}	\|\Psi_{(j)}^n(v^n,\tau)-\Psi_{(j)}^n(v(t_n),\tau)\|_{L^2}.$$
Consequently, this result and Lemmas~\ref{0psi1}--\ref{0psi2-7} yield the desired result.
\end{proof}
\subsubsection{Proof of the result \eqref{T1.2}}
\begin{proof}
	From \eqref{error expression}, we have
	$
	g_{n+1}=g_n+\mathcal L^n+
	\sum_{j\in\{3,5,7\}}
	\big[\Psi_{(j)}^n(v^n,\tau)-\Psi_{(j)}^n(v(t_n),\tau)\big].
	$
	Summing over the time steps and using $g_0=0$ yields
	\[
	g_n=\sum_{m=0}^{n-1}\mathcal L^m+
	\sum_{m=0}^{n-1}\sum_{j\in\{3,5,7\}}
	\big[\Psi_{(j)}^m(v^m,\tau)-\Psi_{(j)}^m(v(t_m),\tau)\big].
	\]
	Combining Lemma~\ref{lem:cumulative-local-error} with the stability
	estimates established in the proof of Lemma~\ref{0stability},
	we obtain
	$
	\|g_n\|_{L^2}\le C\tau^3+C\tau\sum_{m=0}^{n-1}\|g_m\|_{L^2},
	\,1\le n\le L,
	$
	for any $\tau\in(0,\tau_1]$, where $\tau_1$ is as in
	Lemma~\ref{1bound}.
	Applying the discrete Gronwall inequality yields
	\[
	\max_{0\le n\le L}\|v(t_n)-v^n\|_{L^2}\le C_1\tau^3,
	\]
	where $\tau_1$ and $C_1$ depend only on $T$ and
	$\|v\|_{C([0,T];H^{2})}$.
	Since $S_t$ is an isometry on $L^2$, the convergence
	result~\eqref{T1.2} follows immediately.
\end{proof}
\par Having established \eqref{T1.1}--\eqref{T1.2}, the proof of Theorem~\ref{thm1.1} is complete.
\section{Conclusion}\label{section:conclusion}
In this paper, we constructed and analyzed an explicit third-order low-regularity integrator  \texttt{LRI3} for the cubic nonlinear Schr\"odinger equation on the torus $\mathbb{T}$. The proposed scheme utilized a second-order predictor to capture the nonlinear evolution and performed a refined algebraic analysis of the resonant structure to mitigate derivative loss. We established rigorous convergence analysis, proving that \texttt{LRI3} achieves third-order accuracy in $H^r$ for initial data in $H^{r+2}$ with $r>1/2$. Notably, for $H^{2}$ initial data, the method attains third-order convergence in the $L^2$ norm. The numerical simulations presented further substantiated the theoretical convergence rates.


\section*{Funding}
This work was partially supported by the National Natural Science Foundation of China under Grants No. 12371403 (B. W.) and No. 12371395 (Y. W.).

\bibliographystyle{abbrvnat}
\bibliography{reference}

@article{ostermann2022second,
	author    = "Ostermann, Alexander and Wu, Yifei and Yao, Fangyan",
	title="A second-order low-regularity integrator for the nonlinear {S}chr{\"{o}}dinger equation",
	journal   = "Adv. Contin. Discrete Models",
	volume    = "2022",
	number    = "1",
	pages     = "23",
	year      = "2022",
}

@article{cao2024new,
	author		= "Cao, Jiachuan and Li, Buyang and Lin, Yanping",
	title		= "A new second-order low-regularity integrator for the cubic nonlinear {S}chr{\"o}dinger equation",
	journal		= "IMA J. Numer. Anal.",
	volume		= "44",
	number		= "3",
	pages		= "1313--1345",
	year		= "2024",
}

@article{li2021fully,
	author		= "Li, Buyang and Wu, Yifei",
	title		= "A fully discrete low-regularity integrator for the {1D} periodic cubic nonlinear {Schr\"odinger} equation",
	journal		= "Numer. Math.",
	volume		= "149",
	number		= "1",
	pages		= "151--183",
	year		= "2021",
	publisher	= "Springer",
}

@article{ostermann2018low,
	author		= "Ostermann, Alexander and Schratz, Katharina",
	title		= "Low regularity exponential-type integrators for semilinear {S}chr{\"o}dinger equations",
	journal		= "Found. Comput. Math.",
	volume		= "18",
	number		= "3",
	pages		= "731--755",
	year		= "2018",
}

@article{alama2024symmetric,
	author		= "Alama Bronsard, Yvonne",
	title		= "A symmetric low-regularity integrator for the nonlinear {S}chr{\"o}dinger equation",
	journal		= "IMA J. Numer. Anal.",
	volume		= "44",
	number		= "6",
	pages		= "3648--3682",
	year		= "2024",
}

@article{bourgain1993fourier,
	author		= "Bourgain, J.",
	title		= "Fourier transform restriction phenomena for certain lattice subsets and applications to nonlinear evolution equations: {P}art {I}: {S}chr{\"o}dinger equations",
	journal		= "Geom. Funct. Anal.",
	volume		= "3",
	number		= "2",
	pages		= "107--156",
	year		= "1993",
}

@article{kato1988commutator,
	author		= "Kato, Tosio and Ponce, Gustavo",
	title		= "Commutator estimates and the {E}uler and {Navier-Stokes} equations",
	journal		= "Comm. Pure Appl. Math.",
	volume		= "41",
	number		= "7",
	pages		= "891--907",
	year		= "1988",
}

@misc{bronsard2025symmetric,
    author = "Alama Bronsard, Yvonne and Bruned, Yvain and Maierhofer, Georg and Schratz, Katharina",
	title     = "Symmetric resonance based integrators and forest formulae",
	year      = "2025",
	note      = "Preprint at \url{https://arxiv.org/abs/2305.16737}"
}

@article{hochbruck2005explicit,
	author		= "Hochbruck, Marlis and Ostermann, Alexander",
	title			= "Explicit exponential {R}unge--{K}utta methods for semilinear parabolic problems",
	journal		= "SIAM J. Numer. Anal.",
	volume		= "43",
	number		= "3",
	pages			= "1069--1090",
	year			= "2005",
}

@article{bruned2022resonance,
	author    = "Bruned, Yvain and Schratz, Katharina",
	title     = "Resonance-based schemes for dispersive equations via decorated trees",
	journal   = "Forum Math. Pi",
	volume    = "10",
	pages     = "{e2}",
	year      = "2022",
}

@article{Wu2021,
	author		= "Wu, Yifei and Yao, Fangyan",
	title			= "A first-order {Fourier} integrator for the nonlinear {Schr\"odinger} equation on {$\mathbb{T}$} without loss of regularity",
	journal		= "Math. Comp.",
	volume		= "91",
	number="335",
	pages			= "1213--1235",
	year			= "2022",
}

@book{sulem1999,
	author    = "Sulem, C. and Sulem, P.-L.",
	title     = "The Nonlinear {S}chr{\"o}dinger Equation: Self-Focusing and Wave Collapse",
	address   = "New York",
	publisher = "Springer",
	year      = "1999"
}

@book{agrawal2019nonlinear,
	author    = "Agrawal, Govind P.",
	title     = "Nonlinear {F}iber {O}ptics",
	address   = "London",
	publisher = "Academic Press",
	year      = "2019"
}

@article{Ostermann2022,
	author    = "Ostermann, Alexander and Rousset, Frédéric and Schratz, Katharina",
	title     = "Fourier integrator for periodic {NLS}: low regularity estimates via discrete {Bourgain} spaces",
	journal   = "J. Eur. Math. Soc. (JEMS)",
	volume    = "25",
	number    = "10",
	pages     = "3913--3952",
	year      = "2023",
}

@article{ignat2011splitting,
	author    = "Ignat, Liviu I.",
	title     = "A splitting method for the nonlinear {S}chr{\"o}dinger equation",
	journal   = "J. Differential Equations",
	volume    = "250",
	number    = "7",
	pages     = "3022--3046",
	year      = "2011",
}

@article{su2020time,
	author    = "Su, Chunmei and Zhao, Xiaofei",
	title     = "On time-splitting methods for nonlinear {S}chr{\"o}dinger equation with highly oscillatory potential",
	journal   = "ESAIM Math. Model. Numer. Anal.",
	volume    = "54",
	number    = "5",
	pages     = "1491--1508",
	year      = "2020",
}

@article{eilinghoff2016fractional,
	author    = "Eilinghoff, Johannes and Schnaubelt, Roland and Schratz, Katharina",
	title     = "Fractional error estimates of splitting schemes for the nonlinear {S}chr{\"o}dinger equation",
	journal   = "J. Math. Anal. Appl.",
	volume    = "442",
	number    = "2",
	pages     = "740--760",
	year      = "2016",
}

@article{lubich2008splitting,
	author    = "Lubich, Christian",
	title     = "On splitting methods for {S}chr{\"o}dinger--{P}oisson and cubic nonlinear {S}chr{\"o}dinger equations",
	journal   = "Math. Comp.",
	volume    = "77",
	number    = "264",
	pages     = "2141--2153",
	year      = "2008",
}

@article{berland2007conservation,
	author    = "Berland, H{\aa}vard and Islas, Alvaro L. and Schober, Constance M.",
	title     = "Conservation of phase space properties using exponential integrators on the cubic {S}chr{\"o}dinger equation",
	journal   = "J. Comput. Phys.",
	volume    = "225",
	number    = "1",
	pages     = "284--299",
	year      = "2007",
}

@article{celledoni2008symmetric,
	author    = "Celledoni, Elena and Cohen, David and Owren, Brynjulf",
	title     = "Symmetric exponential integrators with an application to the cubic {S}chr{\"o}dinger equation",
	journal   = "Found. Comput. Math.",
	volume    = "8",
	number    = "3",
	pages     = "303--317",
	year      = "2008",
}

@article{hochbruck2010exponential,
	author    = "Hochbruck, Marlis and Ostermann, Alexander",
	title     = "Exponential integrators",
	journal   = "Acta Numer.",
	volume    = "19",
	pages     = "209--286",
	year      = "2010",
}

@article{akrivis1993finite,
	author    = "Akrivis, Georgios D.",
	title     = "Finite difference discretization of the cubic {S}chr{\"o}dinger equation",
	journal   = "IMA J. Numer. Anal.",
	volume    = "13",
	number    = "1",
	pages     = "115--124",
	year      = "1993",
}

@article{henning2017crank,
	author    = "Henning, Patrick and Peterseim, Daniel",
	title     = "{C}rank--{N}icolson {G}alerkin approximations to nonlinear {S}chr{\"o}dinger equations with rough potentials",
	journal   = "Math. Models Methods Appl. Sci.",
	volume    = "27",
	number    = "11",
	pages     = "2147--2184",
	year      = "2017",
}

@article{sanzserna1984methods,
	author    = "Sanz-Serna, J. M.",
	title     = "Methods for the numerical solution of the nonlinear {S}chroedinger equation",
	journal   = "Math. Comp.",
	volume    = "43",
	number    = "167",
	pages     = "21--27",
	year      = "1984",
}

@article{Onorato2001,
	author		= "Onorato, Miguel and Osborne, Alfred R. and Serio, Marina and Bertone, Serena",
	title			= "Freak Waves in Random Oceanic Sea States",
	journal		= "Phys. Rev. Lett.",
	volume		= "86",
	number = "25",
	pages			= "5831--5834",
	year			= "2001",
}

@article{Hofmanov2017,
	author    = "Hofmanová, Martina and Schratz, Katharina",
	title     = "{An exponential-type integrator for the KdV equation}",
	journal   = "Numer. Math.",
	volume    = "136",
	number    = "4",
	pages     = "1117--1137",
	year      = "2017",
}

@article{Ning2022,
	author    = "Ning, Cui and Wu, Yifei and Zhao, Xiaofei",
	title     = "An Embedded Exponential-Type Low-Regularity Integrator for {mKdV} Equation",
	journal   = "SIAM J. Numer. Anal.",
	volume    = "60",
	number    = "3",
	pages     = "999--1025",
	year      = "2022",
}

@article{Ostermann2020,
	author    = "Ostermann, Alexander and Su, Chunmei",
	title     = "{A Lawson-type exponential integrator for the Korteweg–de Vries equation}",
	journal   = "IMA J. Numer. Anal.",
	volume    = "40",
	number    = "4",
	pages     = "2399--2414",
	year      = "2020",
}

@article{Ostermannsu2019,
	author    = "Ostermann, Alexander and Su, Chunmei",
	title     = "Two exponential-type integrators for the “good” {Boussinesq equation}",
	journal   = "Numer. Math.",
	volume    = "143",
	number    = "3",
	pages     = "683--712",
	year      = "2019",
}

@article{Schratz2021,
	author    = "Schratz, Katharina and Wang, Yan and Zhao, Xiaofei",
	title     = "{Low-regularity integrators for nonlinear Dirac equations}",
	journal   = "Math. Comp.",
	volume    = "90",
	number    = "327",
	pages     = "189--214",
	year      = "2021",
}

@article{Wang2022,
	author    = "Wang, Yan and Zhao, Xiaofei",
	title     = "{A symmetric low-regularity integrator for nonlinear Klein-Gordon equation}",
	journal   = "Math. Comp.",
	volume    = "91",
	number    = "337",
	pages     = "2215--2245",
	year      = "2022",
}

@article{Wuzh2022,
	author    = "Wu, Yifei and Zhao, Xiaofei",
	title     = "{Optimal convergence of a second-order low-regularity integrator for the KdV equation}",
	journal   = "IMA J. Numer. Anal.",
	volume    = "42",
	number    = "4",
	pages     = "3499--3528",
	year      = "2022",
}

@article{Wuzh22022,
	author    = "Wu, Yifei and Zhao, Xiaofei",
	title     = "{Embedded exponential-type low-regularity integrators for KdV equation under rough data}",
	journal   = "BIT",
	volume    = "62",
	number    = "3",
	pages     = "1049--1090",
	year      = "2022",
}

@article{Li2023,
	author    = "Li, Hang and Su, Chunmei",
	title     = "A Semi-discrete First-Order Low Regularity Exponential Integrator for the “good” {Boussinesq} Equation Without Loss of Regularity",
	journal   = "J. Sci. Comput.",
	volume    = "95",
	number    = "3",
    pages="74",
	year      = "2023",
}

@article{Ostermannyao2022,
	author    = "Ostermann, Alexander and Yao, Fangyan",
	title     = "A Fully Discrete Low-Regularity Integrator for the Nonlinear {Schr\"odinger} Equation",
	journal   = "J. Sci. Comput.",
	volume    = "91",
	number    = "1",
	pages     = "9",
	year      = "2022",
}

@article{Bao_2024,
  author        = "Bao, Weizhu and Wang, Chushan",
  title         = "An explicit and symmetric exponential wave integrator for the nonlinear {S}chr{\"o}dinger equation with low regularity potential and nonlinearity",
  journal       = "SIAM J. Numer. Anal.",
  volume        = "62",
  number        = "4",
  pages         = "1901--1928",
  year          = "2024",
}

@article{Bao_2024_2,
  author        = "Bao, Weizhu and Wang, Chushan",
  title         = "Optimal error bounds on the exponential wave integrator for the nonlinear {S}chr{\"o}dinger equation with low regularity potential and nonlinearity",
  journal       = "SIAM J. Numer. Anal.",
  volume        = "62",
  number        = "1",
  pages         = "93--118",
  year          = "2024",
}

@article{Zhao_2021,
  author        = "Zhao, Xiaofei",
  title         = "Numerical integrators for continuous disordered nonlinear {S}chr{\"o}dinger equation",
  journal       = "J. Sci. Comput.",
  volume        = "89",
  number        = "2",
  pages="40",
  year          = "2021"
}

\end{document}